\documentclass[3p,preprint]{elsarticle}
\pdfoutput=1

\usepackage{graphicx}
\usepackage{lineno}
\modulolinenumbers[5]

\journal{Journal of Computational Physics}

\usepackage{amsmath, amssymb, amsthm}
\usepackage{booktabs}
\usepackage{nicefrac}
\usepackage{capt-of}
\newcommand{\figurewidthWW}{0.365}    % width of matrix (omega-omega) plots
\newcommand{\figurewidthTR}{0.49}  % width of other plots

\newcommand{\coloneq}{\mathrel{\mathop:}=}
\newcommand{\eqcolon}{=\mathrel{\mathop:}}
\newcommand{\vect}[1]{\mathbf{#1}}
\newcommand{\tvect}[1]{\MakeUppercase{\bf #1}}
\newcommand{\inv}[1]{{#1}^{-1}}
\newcommand{\dt}{\Delta t}
\renewcommand{\d}{\mathrm{d}}
\newcommand{\ii}{\mathrm{i}}
\newcommand{\scrI}{\mathcal{I}}
\newcommand{\scrR}{\mathcal{R}}
\DeclareMathOperator{\spectral}{sp}
\renewcommand{\exp}[1]{\mathrm{e}^{#1}}
\newcommand{\krI}[2]{#1^{#2}}

\biboptions{sort&compress}

\begin{document}

\begin{frontmatter}

\title{A high-order Boris integrator}%\tnoteref{mytitlenote}}

%\author{Author 1}
%\fnref{myfootnote}}
%\address{Address Author 1}

%\author{Author 2}
%\address{Address Author 2}

\author[a2]{Mathias Winkel\corref{cor}}
\ead{mathias.winkel@usi.ch}

\author[a1,a2]{Robert Speck}
\ead{r.speck@fz-juelich.de}

\author[a2]{Daniel Ruprecht}
\ead{daniel.ruprecht@usi.ch}

\address[a2]{Institute of Computational Science, University of Lugano, Switzerland.}
\address[a1]{J\"ulich Supercomputing Centre, Forschungszentrum  J\"ulich, Germany.}

\cortext[cor]{Corresponding author}

\begin{abstract}
This work introduces the high-order Boris-SDC method for integrating the equations of motion for electrically charged particles in an electric and magnetic field.
Boris-SDC relies on a combination of the Boris-integrator with spectral deferred corrections (SDC).
SDC can be considered as preconditioned Picard iteration to compute the stages of a collocation method.
In this interpretation, inverting the preconditioner corresponds to a sweep with a low-order method.
In Boris-SDC, the Boris method, a second-order Lorentz force integrator based on velocity-Verlet, is used as a sweeper/preconditioner.
The presented method provides a generic way to extend the classical Boris integrator, which is widely used in essentially all particle-based plasma physics simulations involving magnetic fields, to a high-order method.
Stability, convergence order and conservation properties of the method are demonstrated for different simulation setups.
Boris-SDC reproduces the expected high order of convergence for a single particle and for the center-of-mass of a particle cloud in a Penning trap and shows good long-term energy stability.
%
%Besides its high formal order, its attractive properties and ease-of-implementation provide the foundation of a considerably improved approach to particle trajectory integration in electric and magnetic fields and pioneers upcoming time-parallel integration methods for this problem.
%Furthermore, the formal methodology presented here can be easily adapted to other specialized integration methods.
\end{abstract}


\end{frontmatter}

%\linenumbers

%%%%%%%%%%%%%%%%%%%%%%%%%%%%%%%%%%%%%%%%%%%%%%%%%%%%%%%%%%%%%%%%%%%%%%%%%%%%%%%%
%%%%%%%%%%%%%%%%%%%%%%%%%%%%%%%%%%%%%%%%%%%%%%%%%%%%%%%%%%%%%%%%%%%%%%%%%%%%%%%%
%%%%%%%%%%%%%%%%%%%%%%%%%%%%%%%%%%%%%%%%%%%%%%%%%%%%%%%%%%%%%%%%%%%%%%%%%%%%%%%%
%%%%%%%%%%%%%%%%%%%%%%%%%%%%%%%%%%%%%%%%%%%%%%%%%%%%%%%%%%%%%%%%%%%%%%%%%%%%%%%%

\section{Introduction}\label{sec:intro}

Often when modeling phenomena in plasma physics, for example particle dynamics in fusion vessels or particle accelerators, an externally applied magnetic field is vital to confine the particles in the physical device~\cite{Stacey2005, Stacey2010}.
In many cases, such as instabilities~\cite{Goldston1995a} and high-intensity laser plasma interaction~\cite{Gibbon2005b}, the magnetic field even governs the microscopic evolution and drives the phenomena to be studied.
Movement of electrically charged particles in an electric and magnetic field is described by the following equations of motion
\begin{subequations}
  \label{eq:base_eq}
  \begin{align}
    \frac{\d\vect{v}}{\d t} &= \vect{f}(\vect{x},\vect{v}) = \alpha\left[\vect{E}(\vect{x},t) + \vect{v}\times\vect{B}(\vect{x},t)\right]\label{eq:rhs_B-field},\\
    \frac{\d\vect{x}}{\d t} &= \vect{v}\label{eq:dxdtv}
\end{align}\label{eq:newton}
\end{subequations}
with the particle position $\vect{x}\in\mathbb{R}^d$, its velocity $\vect{v}\in\mathbb{R}^d$, the magnetic field $\vect{B}(\vect{x},t)\in\mathbb{R}^d$, electric field $\vect{E}(\vect{x},t)\in\mathbb{R}^d$ and the charge to mass ratio $\alpha\in\mathbb{R}$.
  In~\eqref{eq:rhs_B-field}, the well-known Lorentz force $\vect{f}$ depends on $\vect{x}$ and $\vect{v}$ so that a discretizing~\eqref{eq:base_eq} with a standard velocity-Verlet scheme~\cite{Verlet1967,Birdsall1985} reads
\begin{subequations}
	\label{eq:vv}
  \begin{align}
    \vect{x}_{n+1} &= \vect{x}_n + \dt\left(\vect{v}_n + \frac{\dt}{2}\vect{f}(\vect{x}_n,\vect{v}_n)\right)\label{eq:vv_x}\\
    \vect{v}_{n+1} &= \vect{v}_n + \frac{\dt}{2}\left(\vect{f}(\vect{x}_n,\vect{v}_n) + \vect{f}(\vect{x}_{n+1},\vect{v}_{n+1})\right)\label{eq:vv_v}.
  \end{align}
\end{subequations}
with an implicit velocity update step~\eqref{eq:vv_v}.
The Boris integration method~\cite{Boris1970, Birdsall1985} provides a clever way to evaluate~\eqref{eq:vv_v} without having to actually solve an implicit system.
It has thus become a de-facto standard for the numerical solution of~\eqref{eq:vv} and allows to cheaply integrate the particle trajectory in the presence of electric and magnetic fields.
  
Being based on the velocity-Verlet scheme, the Boris approach is a second-order method~\cite{Birdsall1985}.
Whether it is also symplectic is controversial:
In~\cite{Webb2014} it is claimed that it is while~\cite{QinEtAl2013} claim that it is not, but nevertheless shows excellent long term energy stability due to being phase-space volume preserving.
Furthermore, it only requires a single evaluation of the right-hand side $\vect{f}$ per time step, making it a cheap numerical method in terms of computational cost~\cite{Patacchini2009}.
For these reasons, the Boris method is widely used in many Particle-In-Cell-codes (see e.\,g.~\cite{Verboncoeur2005}), grid-free methods (e.\,g.~\cite{Gibbon2010a}) and Monte-Carlo simulations (e.\,g.~\cite{Kirschner2000}).
Several explicit alternatives to the Boris method have been proposed, compare~\cite{Patacchini2009} and references therein.
All of them are second-order accurate and apparently no higher-order methods based on the original Boris approach exist.
 Especially for applications such as trajectory integration in particle accelerators~\cite{Toggweiler2014}, space-weather studies~\cite{Lapenta2012}, high-intensity laser-plasma interaction~\cite{Gibbon2005b}, and fusion vessel simulations~\cite{Bowers2009,Kirschner2000}, where high precision has to be maintained over long physical simulation times, these are desirable, though.
In addition, the current development of high-performance computing systems towards high floating point operation rates at stagnating memory data transfer speed  favors the use of higher-order methods in essentially all fields of computer simulation~\cite{Dongarra2014}. 
Furthermore, the ability of tuning the order of an integration algorithm adds a new dimension to its parameter space that allows for balancing precision versus runtime.
  % as it is for example necessary for parallel-in-time methods, such as Parareal~\cite{MadayTurinici2005} and PFASST~\cite{EmmettMinion2012}.
  
A number of other high-precision or higher-order methods for~\eqref{eq:base_eq} have been developed.
Examples are methods that use a spatial coordinate instead of time as the independent variable which showed better performance than a fourth-order Runge-Kutta method in beam propagation simulations~\cite{Stoltz2002}, a quasi-symplectic Trotter-factorization based scheme that builds upon an explicit-implicit mixture of leap-frog, Verlet, exponential differencing and Boris rotations with a sixth-order rotation angle approximation~\cite{Bowers2009} or a Taylor series-based explicit approach with an up to sixth-order replacement for the Boris method using a complex differential operator for the Maxwell fields~\cite{Quandt2010, Quandt2007}.
Essentially, none of these methods are easily tunable for arbitrary order but are formulated for a very specific case.
Only the latter, Taylor series-based approach offers this feature but requires a complicated set of appropriate differential operators to be constructed for every order.

%  WHAT WE ARE OFFERING
In this work, we introduce the high-order \emph{Boris-SDC} integration method for~\eqref{eq:newton}, which is a combination of the classical Boris integrator with the spectral deferred corrections method (SDC).
The resulting Boris-SDC method retains much of the simplicity of the Boris integrator (in terms of implementation, alas not derivation) while allowing to easily generate a method of essentially arbitrary order.
Based on classical defect correction, SDC has originally been introduced for first-order ODEs as an iterative approach for the generic construction of high-order integration schemes using a low-order base propagator (the ``sweeper'') such as implicit or explicit Euler for the correction ``sweeps''~\cite{DuttEtAl2000}.
Several modifications and extensions exist e.\,g.~semi-implicit SDC~\cite{Minion2003}, GMRES-accelerated SDC~\cite{HuangEtAl2006}, inexact SDC~\cite{SpeckEtAl2014_DDM2013}, multi-level SDC~\cite{SpeckEtAl2014_BIT}, SDC based on DIRK methods~\cite{Weiser2014} or the parallel full approximation scheme in space and time (PFASST), a parallel-in-time integrator which exploits the iterative structure of SDC~\cite{EmmettMinion2012}.
Recently, SDC has been formulated for second-order problems with the standard Verlet integrator as sweeper~\cite{Minion2014_2ndOrderSDC}.
Here, we combine this particular approach with the classical Boris integrator and extend it to a velocity-dependent force of the form~\eqref{eq:rhs_B-field}.

The derivation of Boris-SDC relies on the interpretation of SDC as a preconditioned Picard iteration for the solution of a collocation problem, see e.\,g.~\cite{HuangEtAl2006,MinionEtAl2015}.
Collocation methods are based on the integral formulation of an ODE, the approximation of the exact trajectory over a time step by a polynomial and evaluation of the integrals by quadrature.
They are a special class of implicit Runge-Kutta methods and, depending on the chosen quadrature nodes, have a number of attractive properties, particularly symplecticness, see e.\,g.~\cite{hairer_nonstiff,hairer_stiff}.
The disadvantage of collocation methods is that they require the solution of a very large, possibly nonlinear system of equations to compute the stages.  
Picard iteration can be used to solve this system, but often requires a too small time step for convergence.
SDC can be considered as a preconditioned Picard iteration, where inverting the preconditioner corresponds to ''sweeping'' through the quadrature nodes with a low-order method.
If sufficiently many sweeps are performed, the advantageous properties of the underlying collocation method are recovered.
For e.\,g.~a first-order method such as the implicit Euler as sweeper, SDC formally gains one order per sweep~\cite{ShuEtAl2007}, so fixing the number of iterations allows to easily generate a scheme of higher order, up to the order provided by the underlying quadrature.
Here, we describe how the classical Boris integrator can be used as a preconditioner to derive an iterative solver for a collocation approximation of~\eqref{eq:base_eq}.
% und this a high-order integrator for~\eqref{eq:base_eq}.

% PREVIEW ON THIS WORK
This paper is organized as follows:
Section~\ref{sec:coll} describes collocation methods, briefly discusses their properties and introduces the required notation.
In Section~\ref{sec:borissdc}, we start with spectral deferred corrections based on the velocity-Verlet scheme as base integrator in matrix from.
The matrix formulation itself is derived in~\ref{sec:apx_vv}.
Concentrating on this rather formal notation, these parts are sufficiently general to also be utilized for force expressions other than the Lorentz force in~\eqref{eq:rhs_B-field}.
In the second part of Section~\ref{sec:borissdc} we then specialize the formalism to the case of the Lorentz force as the ODE's right-hand side and derive ready-to-implement expressions for the Boris-SDC method, specifically tailored for problems of the form~\eqref{eq:base_eq}.
Section~\ref{sec:numerics} illustrates the properties of Boris-SDC by numerical examples and compares Boris-SDC to the classical Boris integrator. Finally, Section~\ref{sec:conclusion} gives a summary and an outlook on possible future directions of research.

\section{Collocation formulation}\label{sec:coll}

In this section, we briefly describe collocation methods and introduce the notation required for the spectral deferred correction approach of Section~\ref{sec:borissdc}.
Note that the notation below is loosely based on the discussion of collocation methods and SDC with a velocity-Verlet integrator as base method for second-order problems in~\cite{Minion2014_2ndOrderSDC}.

Rewriting equations~\eqref{eq:base_eq} in Picard formulation for an arbitrary interval $[t_n,t_{n+1}]$ with starting value $\vect{x}_0 = \vect{x}(t_n)$, $\vect{v}_0 = \vect{v}(t_n)$, $\vect{x},\vect{v}\in\mathbb{R}^d$, we obtain%
\footnote{We use the following notation here: 
Vectors and matrices in normal font refer to scalar values at a single node in time.
Bold-faced variables indicate aggregation over spatial variables (e.\,g.~particles in more than one dimensions and/or multiple particles).
Vectors with capitals are used to denote aggregation over all intermediate steps.
Matrices are always denoted with capital letters, slightly abusing our own convention.
All matrices in this context, however, refer to aggregated quantities anyway.}
\begin{subequations}
\label{eq:picard}
\begin{align}
    \vect{v}(t) &= \vect{v}_0 + \int_{t_n}^t \vect{f}(\vect{x}(s),\vect{v}(s))\ \d s, \label{eq:v} \\
    \vect{x}(t) &= \vect{x}_0 + \int_{t_n}^t \vect{v}(s)\ \d s. \label{eq:x}
\end{align}
\end{subequations}
Collocation methods are based on the introduction of intermediate nodes
\begin{align}
t_n  \le \tau_1<\ldots<\tau_M\le t_{n+1}, \ M\ge 1
\end{align}
and approximating the integrals using quadrature. 
Details can be found e.\,g.~in~\cite[II.7]{hairer_nonstiff}. 
To allow for a more convenient notation below, we set $\tau_0 \coloneq t_n$.
The quadrature weights are collected in a matrix $\bar{Q}\in\mathbb{R}^{M \times M}$ with entries
\begin{align}
	\bar{q}_{m,j} \coloneq \int_{t_n}^{\tau_m} \ell_j(s) \ \d s, \quad m,j = 1, \ldots, M
\end{align}
where $\ell_j(s)$, $j= 1,\ldots, M$ are Lagrange polynomials.
Then, in order to account for the initial values $\vect{x}_0$, $\vect{v}_0$, we further define
\begin{align}
	Q \coloneq \begin{pmatrix}
		0 & \tvect{0} \\
		\tvect{0} & \bar{Q}
	\end{pmatrix} 
	\in \mathbb{R}^{(M+1)\times(M+1)}.
\end{align}
For quadrature nodes with $\tau_0 = \tau_1$, e.\,g.~Gauss-Lobatto nodes, the second row of $Q$ is zero as well, because $t_n = \tau_0 = \tau_1$.

Let $\tvect{v} = \left(\vect{v}_0,\vect{v}_1,\ldots,\vect{v}_M\right)^T,\ \tvect{x} = \left(\vect{x}_0,\vect{x}_1,\ldots,\vect{x}_M\right)^T\in\mathbb{R}^{(M+1)d}$ be vectors with approximate values for $\vect{x}$ and $\vect{v}$ at the nodes $\tau_m$ and $\tvect{f}(\tvect{x},\tvect{v}) = \left(\vect{f}_0,\vect{f}_1,,\ldots,\vect{f}_M\right)^T\in\mathbb{R}^{(M+1)d}$ with $\vect{f}_m = \vect{f}(\vect{x}_m, \vect{v}_m)$ the vector containing the corresponding right-hand side values.
Then, the matrix $\vect{Q} \coloneq Q\otimes\vect{I}_d$ provides approximations of the integral in~\eqref{eq:v} with $t=\tau_m$, that is
\begin{align}
	\label{eq:disc_int}
    \vect{Q}\tvect{F} (\tvect{x},\tvect{v}) =(Q\otimes\vect{I}_d)\tvect{f}(\tvect{x},\tvect{v}) \approx \left(\int_{t_n}^{\tau_m} \vect{f}(\vect{x}(s),\vect{v}(s))\ \d s\right)_{m=0,1,\ldots,M}
\end{align}
with $\vect{I}_d\in\mathbb{R}^{d\times d}$ being the identity matrix and $\otimes$ the standard Kronecker product.
Similarly, the term $\vect{Q}\tvect{v}$ approximates the integrals over $\vect{v}$.
The discrete version of ~\eqref{eq:picard} is then given by the collocation formulation
\begin{subequations}
\label{eq:deq}
\begin{align}
    \tvect{v} &= \tvect{v}_0 + \vect{Q}\tvect{f}(\tvect{x},\tvect{v}) \label{eq:veldeq}\\
    \tvect{x} &= \tvect{x}_0 + \vect{Q}\tvect{v} = \tvect{x}_0 + \vect{Q}\tvect{v}_0 + \vect{Q}\vect{Q}\tvect{f}(\tvect{x},\tvect{v}). \label{eq:spcdeq}
\end{align}
\end{subequations}
for 
\begin{align}
\tvect{x}_0 \coloneq \left(\vect{x}_0,\vect{x}_0,\ldots,\vect{x}_0\right)^T, \tvect{v}_0 \coloneq \left(\vect{v}_0,\vect{v}_0,\ldots,\vect{v}_0\right)^T\in\mathbb{R}^{(M+1)d}.
\end{align}
We note that $\vect{QQ} = (Q\otimes\vect{I}_d)(Q\otimes\vect{I}_d) = QQ\otimes \vect{I}_d$.

\subsection{Collocation system}

In order to combine both equations into a single, closed expression, we first note that $\tvect{f}$ depends on both $\tvect{x}$ and $\tvect{v}$, or, more precisely, each component $\vect{f}_m$ depends on the tuple $(\vect{x}_m,\vect{v}_m)\in\mathbb{R}^{2d}$.
The ordering of Equations~\eqref{eq:deq} (first all the $\vect{v}_m$, then all the $\vect{x}_m$), however, is not compatible with the sorting in $\tvect{f}$, where the first entry is depends on $(\vect{x}_1, \vect{v}_1)$, the second on $(\vect{x}_2, \vect{v}_2)$ etc.
Thus, we need to resort the matrix formulation~\eqref{eq:deq} so that the degrees-of-freedom are ordered as 
\begin{align}
    \label{eq:apx_mixed_vec}
    \tvect{u} = (\vect{u}_0,\vect{u}_1,\ldots,\vect{u}_M)^T \coloneq (\vect{x}_0,\vect{v}_0,\ldots,\vect{x}_m,\vect{v}_m,\ldots,\vect{x}_M,\vect{v}_M)^T\in\mathbb{R}^{(M+1)2d}.
\end{align}
To this end, we introduce permutation operators $\krI{I}{x}$, $\krI{I}{v}$ and $\krI{I}{xv}$ with
\begin{align}
 	\krI{I}{x} = 
 	\begin{pmatrix}
 	1 \\ 0
 	\end{pmatrix},\quad \krI{I}{v} = 
 	\begin{pmatrix}
 	0 \\ 1
 	\end{pmatrix}\quad\text{and}\quad \krI{I}{xv} =
 	\begin{pmatrix}
 	0 & 1\\
 	0 & 0
 	\end{pmatrix}\label{eq:I_matrices}
\end{align}
so that $\vect{u}_m = \left( \vect{I}_{d} \otimes \krI{I}{x} \right) \vect{x}_m + \left( \vect{I}_{d} \otimes \krI{I}{v} \right) \vect{v}_m$ for $m=0,\ldots,M$. 
The permutation operators $\krI{I}{x}$ and $\krI{I}{v}$ redistribute the entries of $\vect{x}_{m}$ and $\vect{v}_{m}$ to match the sorting of the degrees of freedom in $\tvect{u}$ while the operator $\krI{I}{xv}$ reflects the action of a velocity component (entry in the second column) on a position component (first row).

For a matrix $R\in\mathbb{R}^{n\times n}$ we now use for abbreviation $\krI{R}{x} \coloneq R\otimes \krI{I}{x}\in\mathbb{R}^{2n\times n}$, $\krI{R}{v} \coloneq R\otimes \krI{I}{v}\in\mathbb{R}^{2n\times n}$, and $\krI{R}{xv} \coloneq R\otimes \krI{I}{xv}\in\mathbb{R}^{2n\times 2n}$.
Then, the redistributed version of~\eqref{eq:deq} reads
\begin{align}\label{eq:u_def_coll}
	\tvect{u} = \tvect{u}_0 + \vect{\krI{Q}{xv}}\tvect{u}_0 + \vect{\krI{QQ}{x}}\tvect{f}(\tvect{U}) + \vect{\krI{Q}{v}}\tvect{f}(\tvect{U})
\end{align}
with $\tvect{u}_0 = (\vect{u}_0,\ldots,\vect{u}_0)^T\in\mathbb{R}^{(M+1)2d}$, $\tvect{f}(\tvect{U}) \coloneq \tvect{f}(\tvect{X},\tvect{V})$
and
\begin{subequations}
\begin{align}
	\vect{\krI{Q}{v}} &\coloneq \krI{Q}{v}\otimes\vect{I}_d = Q\otimes \krI{I}{v}\otimes\vect{I}_d \\
	\vect{\krI{QQ}{x}} &\coloneq QQ\otimes \krI{I}{x}\otimes\vect{I}_d
\end{align}
\end{subequations}

Equation~\eqref{eq:u_def_coll} can be compactly written as a possibly non-linear system of equations
\begin{align}
   \label{eq:pic_system}
    \vect{M}_\mathrm{coll}(\tvect{u}) = \vect{C}_\mathrm{coll}\tvect{u}_0.
\end{align}
with 
\begin{align}
	\vect{C}_\mathrm{coll} \coloneq \vect{I}_{(M+1)2d} +\vect{\krI{Q}{xv}},\quad
	\vect{Q}_\mathrm{coll} \coloneq \vect{\krI{QQ}{x}} + \vect{\krI{Q}{v}},\quad
	\vect{M}_\mathrm{coll}(\cdot) \coloneq \left(\vect{I}_{(M+1)2d} - \vect{Q}_\mathrm{coll}\tvect{f}\right)(\cdot).
\end{align}
Here, $\vect{C}_\mathrm{coll}$ and $\vect{Q}_\mathrm{coll}$ are matrices, while $\vect{M}_\mathrm{coll}$ can in general be a non-linear operator.
Setting
\begin{align}
  \vect{P}_\mathrm{coll} &\coloneq \inv{\vect{M}}_\mathrm{coll}\vect{C}_\mathrm{coll}\label{eq:pcoll},
\end{align}
the formal update for $\tvect{u}$ simply reads $\tvect{u} = \vect{P}_\mathrm{coll}\tvect{u}_0$.
We note that this formalism can be easily extended for higher-order ODEs:
For an $L$th-order ODE formulated as first-order system (as done here for $L=2$), the permutation operators~\eqref{eq:I_matrices} are simply the unit vectors and $\krI{I}{xv}$ is replaced by at set of matrices that couple the components accordingly.

Evaluating $\vect{P}_\mathrm{coll}$ requires the inversion of $\vect{M}_\mathrm{coll}$. 
Only for the sake of notational simplicity, we now focus on linear right-hand side functions $\vect{f}$, so that $\vect{M}_\mathrm{coll}$ is a matrix with
\begin{align}
    \vect{M}_\mathrm{coll}(\cdot) = \vect{M}_\mathrm{coll} = \vect{I}_{(M+1)2d} - \vect{Q}_\mathrm{coll}\tvect{f}
\end{align}
and inverse $\inv{\vect{M}}_\mathrm{coll}$.
However, we emphasize that the very same ideas and formulas described in the following apply for the case of non-linear functions, too, as e.\,g.~shown in Section~\ref{sec:numerics}.
Then, operators like $\inv{\vect{M}}_\mathrm{coll}$ have to be interpreted accordingly, c.\,f.~the discussion at the end of~\ref{sec:apx_vv}.
For clarity, arguments of the (now linear) mapping $\tvect{f}$ are still shown with brackets.

In order to obtain a closed update formula which directly maps the initial data $\vect{u}_0$ to the final value $\vect{u}_{n+1} = (\vect{x}_{n+1},\vect{v}_{n+1})$, we define the linear transfer operators $\vect{T}_\mathrm{P}\in\mathbb{R}^{(M+1)2d\times 2d}$ and $\vect{T}_\mathrm{R}\in\mathbb{R}^{2d\times (M+1)2d}$ via
\begin{align}
	\vect{T}_\mathrm{P}\vect{u}_0 = \tvect{u}_0\quad\ \text{and}\quad\ \vect{T}_\mathrm{R}\tvect{u} = \vect{u}_M.\label{eq:transfer_ops} 
\end{align}
Since $\tau_M$ is not necessarily the final step (if $\tau_M<t_{n+1}$ as e.\,g.~for Gauss-Legendre nodes), we make use of the collocation formulation again to obtain approximations $\vect{v}_{n+1}, \vect{x}_{n+1}$ to the final values $\vect{v}(t_{n+1}),\vect{x}(t_{n+1})$ from the full vector $\tvect{u}$ or $\tvect{f}$, respectively.
We define by $q \coloneq \left(0, \bar{q} \right) = \left(0, \bar{q}_{1}, \ldots, \bar{q}_{M} \right) \in\mathbb{R}^{1\times(M+1)}$ the extended vector of quadrature weights over the full interval $[t_n,t_{n+1}]$, where
\begin{align}
    \bar{q}_m = \int_{t_n}^{t_{n+1}}\ell_m(\tau)\ \d\tau.
\end{align}
Then we obtain
\begin{align}
    \vect{v}_{n+1} &= \vect{v}_0 + \vect{q}\tvect{f}(\tvect{U}),\\
    \vect{x}_{n+1} &= \vect{x}_0 +\vect{q}\tvect{v} = \vect{x}_0 +\vect{q}\tvect{v}_0 + \vect{q}\vect{Q}\tvect{f}(\tvect{U}),
\end{align}
which can be combined into a single equation again using
\begin{align}
    \vect{u}_{n+1} = \tilde{\vect{C}}_\mathrm{coll}\tvect{u}_0 + \tilde{\vect{Q}}_\mathrm{coll}\tvect{f}(\tvect{U})\label{eq:final_step_coll}
\end{align}
with
\begin{align}\label{eq:tilde_defs}
    \tilde{\vect{C}}_\mathrm{coll} \coloneq \vect{T}_\mathrm{R} + \vect{\krI{q}{xv}}\quad\ \text{and}\quad\ \tilde{\vect{Q}}_\mathrm{coll} \coloneq \vect{\krI{qQ}{x}} + \vect{\krI{q}{v}}.
\end{align}
Note that if $\tau_M=t_{n+1}$, the vector $\vect{q}$ is equal to the last row of the matrix $\vect{Q}$ and computing $\vect{x}_{n+1},\vect{v}_{n+1}$ is equivalent to computing $\vect{x}_M,\vect{v}_M$.
Now, the complete update formula for $\vect{u}$ reads
\begin{align}
	\label{eq:coll_update}
    \vect{u}_{n+1} = \tilde{\vect{P}}_\mathrm{coll}(\vect{u}_0) \coloneq \tilde{\vect{C}}_\mathrm{coll}\vect{T}_\mathrm{P}\vect{u}_0 + \tilde{\vect{Q}}_\mathrm{coll}\tvect{f}\left(\vect{P}_\mathrm{coll}\vect{T}_\mathrm{P}\vect{u}_0\right).
\end{align} 
The subsequent parts of this section deal with the numerical properties of this formulation and point towards strategies for efficiently inverting $\vect{M}_\mathrm{coll}$, i.\,e.~solving~\eqref{eq:pic_system} by an iterative method.

\subsection{Properties of collocation methods}\label{sec:collprop}

A collocation method with $M$ nodes is equivalent to an $M$-stage implicit Runge Kutta method (IRKM) with a Butcher tableau
\begin{align}
\begin{tabular}{ c | c  }
$c$ & $ \bar{Q}$   \\ \hline
& $\bar{q} $ \\
\end{tabular}
\end{align}
with $c$ being the vector of nodes $\tau_m$ scaled to the unit interval, see e.\,g.~\cite[Theorem 7.7]{hairer_nonstiff}. 
In this interpretation, equation~\eqref{eq:pic_system} is a system of equations to be solved for the $M$ stages of an IRKM while~\eqref{eq:final_step_coll} is the actual update step to be performed once the stages are known.

Collocation methods have a number of attractive numerical properties:
They are of optimal order, $2 M$ for Legendre and $2 M-2$ for Lobatto nodes.
For both Gauss-Legendre and Gauss-Lobatto nodes, the resulting method is symmetric because the corresponding nodes are symmetric~\cite[Theorem 8.9]{hairer_geometric}.
Also, for Gauss-Legendre nodes, the resulting method is always symplectic~\cite[Theorem 16.5]{hairer_nonstiff} as well as B- and A-stable~\cite[Theorem 12.9]{hairer_stiff}.
In Section~\ref{sec:numerics}, we show that collocation methods with Lobatto nodes also have excellent stability properties for the Penning trap example considered there. 
For Lobatto nodes, however, the method is not necessarily symplectic~\cite[Table 16.2]{hairer_nonstiff}.
Nevertheless, for the cases studied here, the Hamiltonian is given as a quadratic form with a symmetric real matrix for which symmetric methods are also symplectic and vice versa~\cite[Theorem 4.9]{hairer_geometric}. 
Hence, for the problems in Section~\ref{sec:numerics}, Gauss-Lobatto nodes also yield a symplectic collocation method and because they do not require an additional step to compute the final value at the end of the interval, we focus on Lobatto nodes here.
For other cases, Legendre nodes might have to be used to obtain a symplectic method.
Note that despite the collocation method being symplectic, energy drift can still emerge due to accumulation of round-off errors, see~\cite{Hairer2008}.

\subsection{Picard iteration and preconditioning}
%\todoRS{Define residual $r$ (and refer to this from Section~\ref{sec:residual_control})}
Even in the linear case, the dense structure of $\tvect{M}_{\mathrm{coll}}$ calls for an an iterative approach to solve~\eqref{eq:pic_system}.
The simplest iteration procedure is a Richardson iteration, see e.\,g.~\cite{Kelley1995}, reading
\begin{subequations}
\label{eq:discrete_picard}
\begin{align}
    \tvect{u}^0 &= \tvect{u}_0\\
	\tvect{u}^{k+1} &= \left( \vect{I}_{(M+1)2d} - \vect{M}_\mathrm{coll} \right)\tvect{u}^k + \vect{C}_\mathrm{coll} \tvect{u}_0 = \vect{C}_\mathrm{coll}\tvect{u}_0 + \vect{Q}_\mathrm{coll}\tvect{f}(\tvect{u}^k). 
\end{align}
\end{subequations}
for $k=0,\ldots,K$.
Here, superscript $k$ denotes the iteration steps.
As a measure for convergence, the norm $r=||\tvect{r}^{k}||$ of the residual 
\begin{equation}
	\label{eq:residual}
	\tvect{r}^{k}:= \vect{C}_{\mathrm{coll}} \tvect{u}_{0} - \vect{M}_\mathrm{coll} \tvect{u}^{k}
\end{equation}
can be monitored, where $\tvect{r}^k\in\mathbb{R}^{M2d}$.
We note that~\eqref{eq:discrete_picard} is equivalent to a discretized Picard iteration. 
Convergence depends on the eigenvalues of the iteration matrix $\vect{K}_\mathrm{pic} \coloneq \vect{Q}_\mathrm{coll}\tvect{f}$.
In the non-linear case, i.\,e.~where $\vect{K}_\mathrm{pic}$ is an operator, convergence properties are given by some adequate norm of this operator.
%the iteration operator of the Picard iteration. 
Picard iteration typically converges only for very small time steps and is thus usually not an efficient approach, so that more advanced methods based on preconditioners are necessary.
For a concise overview of iterative methods including the concept of preconditioning we refer to~\cite{Kelley1995}.
%The Picard iteration  is directly linked to a splitting of the system operator $\vect{M}_\mathrm{coll}$, where
%\begin{align}
%	\vect{M}_\mathrm{coll} = \vect{M}_\mathrm{coll}-\vect{I}_{(M+1)2d} + \vect{I}_{(M+1)2d}.
%\end{align}
%The first part, $\vect{M}_\mathrm{coll}-\vect{I}_{(M+1)2d}$, is the part of $\vect{M}_\mathrm{coll}$ which is ``hard'' to invert, while the second part, $\vect{I}_{(M+1)2d}$, is the ``easy'' component.
%The iteration operator is then the negative of the first part and the second part is multiplied to the left-hand side.
%Thus, a straightforward extension of this iterative scheme is a preconditioned iteration, where $\vect{M}_\mathrm{coll}$ is split differently.
%We introduce a splitting 
%\begin{align}
%    \vect{M}_\mathrm{coll} = \vect{M}_\mathrm{coll}-\vect{M}_{\mathrm{pc}} + \vect{M}_{\mathrm{pc}},
%\end{align}
Introduction of a preconditioner $\vect{M}_{\mathrm{pc}} \approx \vect{M}_{\mathrm{coll}}$ leads to the preconditioned Richardson iteration
\begin{equation}
	\label{eq:precond_iteration}
	\vect{M}_{\mathrm{pc}} \tvect{u}^{k+1} = \left( \vect{M}_{\mathrm{pc}} - \vect{M}_{\mathrm{coll}} \right) \tvect{u}^{k} + \vect{C}_{\mathrm{coll}} \tvect{u}_{0}.
\end{equation}
Here, each iteration requires to solve a linear or nonlinear system of equations determined by the preconditioner $\vect{M}_{\mathrm{pc}}$.
The key is to find a good  preconditioner:
It has to be easy to invert so that computing~\eqref{eq:precond_iteration} is cheap but still provides a sufficiently good approximation of $\vect{M}_{\mathrm{coll}}$ to lead to robust and rapid convergence.
In the next section, we will construct such a preconditioner out of the well-known Boris integrator for problems of the form~\eqref{eq:base_eq}.

\section{Spectral deferred corrections based on the Boris-integrator (Boris-SDC)}\label{sec:borissdc}
The idea to interpret spectral deferred corrections as a preconditioned iterative scheme has been used for different purposes e.\,g.~in~\cite{HuangEtAl2006,Weiser2014}.
Here, we employ it to derive a problem-specific formulation of SDC based on the Boris integration method.
The formulation follows the derivation of SDC for second-order problems in~\cite{Minion2014_2ndOrderSDC}, using the standard velocity-Verlet integrator as preconditioner.
There, problems are considered with a force field $\vect{f}$ that depends only on the position $\vect{x}$.
This corresponds to a separable Hamiltonian of a specific form.
In general, an $\vect{f}$ that also depends on the velocity $\vect{v}$ leads to an implicit update for $\vect{v}$ in the velocity-Verlet integrator, cf.~\eqref{eq:vv}.
For the specific form of $\vect{f}$ in~\eqref{eq:base_eq}, however, the Boris integrator provides a trick that essentially allows for a very efficient solution of the implicit system.

In order to use the Boris integrator as a preconditioner, we need a formulation of the velocity-Verlet scheme in matrix form, similar to~\eqref{eq:pic_system}.
A concise summary of this rather tedious derivation is given in~\ref{sec:apx_vv}. 

\subsection{Velocity-Verlet-based spectral deferred corrections}\label{sec:vv_sdc}

The standard velocity-Verlet integrator~\eqref{eq:vv} for time steps $\tau_0,\tau_1,\ldots,\tau_M$, $M\ge1$, can be written as system of equations
\begin{align}
	\tvect{M}_{\mathrm{vv}} \tvect{u} = \vect{C}_\mathrm{vv}\tvect{u}_0 \label{eq:vv_as_matrix},
\end{align}
with system matrix
\begin{align}
	\vect{M}_\mathrm{vv} \coloneq \vect{I}_{(M+1)2d} -  \vect{Q}_\mathrm{vv} \tvect{f},
\end{align}
see~\ref{sec:apx_vv}, in particular~\eqref{eq:apx_vv_system} and~\eqref{eq:apx_vv_systemoperator}, for details (we assume linear right-hand side here as well for notational simplicity). 
Since both $\vect{I}_{(M+1)2d}$ and $\vect{Q}_\mathrm{vv}$ are lower block-diagonal matrices,~\eqref{eq:vv_as_matrix} can easily be solved by forward substitution.
Note that for non-linear functions $\vect{f}$, $\vect{M}_\mathrm{vv} = \vect{M}_\mathrm{vv}(\cdot)$ can no longer be compactly written in matrix form, but the system itself is still straightforward to solve.
Hence, $\vect{M}_\mathrm{vv}$ satisfies the conditions necessary for a suitable preconditioner:
It  can be easily inverted and approximates the original action of $\vect{M}_\mathrm{coll}$, as both systems correspond to equations for approximations of the values of $\vect{x}$ and $\vect{v}$ at the quadrature nodes, a high-order approximation in case of~\eqref{eq:pic_system} and a composite low-order approximation in case of~\eqref{eq:vv_as_matrix}.

In order to precondition the iteration~\eqref{eq:discrete_picard}, we therefore apply the splitting
\begin{align}\label{eq:M_split}
\vect{M}_\mathrm{coll} = \vect{M}_\mathrm{coll}-\vect{M}_\mathrm{vv}+\vect{M}_\mathrm{vv}.
\end{align}
Using $\vect{M}_{\mathrm{pc}} := \vect{M}_{\mathrm{vv}} \approx \vect{M}_{\mathrm{coll}}$ as a preconditioner in~\eqref{eq:precond_iteration} results in the iteration
\begin{subequations}
\label{eq:full_matrix_sdc}
\begin{align}
    \tvect{u}^{0} &= \tvect{u}_0\\
	(\vect{I}_{(M+1)2d} - \vect{Q}_\mathrm{vv}\tvect{f}) (\tvect{u}^{k+1}) &= \left(\vect{Q}_\mathrm{coll}-\vect{Q}_\mathrm{vv}\right)\tvect{f}(\tvect{u}^k) + \vect{C}_\mathrm{coll}\tvect{u}_0
\end{align}   
\end{subequations}
for $k=0,\ldots,K$ with iteration matrix
\begin{align}
\vect{K}_\mathrm{sdc} \coloneq  \inv{\left(\vect{I}_{(M+1)2d} -  \vect{Q}_\mathrm{vv}\tvect{f}\right)}\left(\vect{Q}_\mathrm{coll}-\vect{Q}_\mathrm{vv}\right)\tvect{f}.
\end{align}
As the name suggests, this particular choice of preconditioner leads to the method of spectral deferred corrections (SDC) with velocity-Verlet as base integrator.
The ``direct'' update matrix $\vect{P}_\mathrm{coll}$ of~\eqref{eq:pcoll} is then replaced by the sum of the updates given by the preconditioned iteration, i.\,e.~we have
\begin{align}
\label{eq:sdc_update}
\vect{P}^k_\mathrm{sdc} &\coloneq \vect{K}^k_\mathrm{sdc} + \sum_{l=0}^{k-1}\vect{K}_\mathrm{sdc}^l\inv{\vect{M}}_\mathrm{vv}\vect{C}_\mathrm{coll}
\end{align}
so that an approximation $\vect{u}^k_{n+1}$ to $\vect{u}_{n+1}$ with $k$ SDC iterations can be computed through 
\begin{align}
    \vect{u}^k_{n+1} = \tilde{\vect{P}}^k_\mathrm{sdc}\vect{u}_0 \coloneq \tilde{\vect{C}}_\mathrm{coll}\vect{T}_\mathrm{P}\vect{u}_0 + \tilde{\vect{Q}}_\mathrm{coll}\tvect{f}\left(\vect{P}^k_\mathrm{sdc}\vect{T}_\mathrm{P}\vect{u}_0\right).\label{eq:psdc}
\end{align}
Convergence of the SDC iteration can be monitored via the residual~\eqref{eq:residual}.
%On the other hand, with velocity-Verlet as base method, each iteration formally adds two orders of accuracy to the resulting method, up to the upper bound defined by the underlying collocation method (e.\,g.~$2M-2$ for Gauss-Lobatto collocation nodes).

While the formulation with iteration and update matrices is convenient to analyze, a different approach is needed for an actual implementation to avoid the explicit use and storage of the full right-hand side vector $\tvect{f}(\tvect{u})$.
To this end, $\vect{M}_\mathrm{coll}$ is again split according to~\eqref{eq:M_split}, but the components $\tvect{x}$ and $\tvect{v}$ are then treated separately, so that~\eqref{eq:full_matrix_sdc} becomes
\begin{subequations}
\begin{align}
    \tvect{x}^{k+1} -  \vect{Q}_x\tvect{f}(\tvect{u}^{k+1}) &= \left(\vect{Q}\vect{Q}-\vect{Q}_x\right)\tvect{f}(\tvect{u}^k) + \vect{Q}\tvect{v}_0 + \tvect{x}_0,\\
    \tvect{v}^{k+1} -  \vect{Q}_T\tvect{f}(\tvect{u}^{k+1}) &= \left(\vect{Q}-\vect{Q}_T\right)\tvect{f}(\tvect{u}^k) + \tvect{v}_0,
\end{align}
\end{subequations}
using the definitions of $\vect{Q}_T$ and $\vect{Q}_x$ of~\eqref{eq:QT} and~\eqref{eq:Qx} in~\ref{sec:apx_vv}.
For a component-wise notation (in terms of nodes $\tau_1,\ldots,\tau_M$), we define
\begin{align}
    \label{eq:Q_matrices}
    \vect{Q}         &\eqcolon (q_{m,l})_{m,l=0,\ldots,M},   &
    \vect{Q}\vect{Q} &\eqcolon (qq_{m,l})_{m,l=0,\ldots,M}, &
    \vect{Q}_x       &\eqcolon (q^x_{m,l})_{m,l=0,\ldots,M}, &
    \vect{Q}_T       &\eqcolon (q^T_{m,l})_{m,l=0,\ldots,M}    
\end{align}
and obtain
\begin{subequations}
\label{eq:node_to_zero_q}
\begin{alignat}{6}
\label{eq:node_to_zero_v}
\vect{v}^{k+1}_{m+1}  &= \vect{v}_0                                   &&+ \sum_{l=0}^{M}q^T_{m+1,l}\left(\vect{f}(\vect{x}^{k+1}_l,\vect{v}^{k+1}_l) - \vect{f}(\vect{x}^{k}_l,\vect{v}^{k}_l)\right) &&+ \sum_{l=0}^Mq_{m+1,l}\vect{f}(\vect{x}^{k}_l,\vect{v}^{k}_l),\\
\label{eq:node_to_zero_x}
\vect{x}^{k+1}_{m+1}  &= \vect{x}_0 + \sum_{l=0}^Mq_{m+1,l}\vect{v}_0 &&+ \sum_{l=0}^{M}q^x_{m+1,l}\left(\vect{f}(\vect{x}^{k+1}_l,\vect{v}^{k+1}_l) - \vect{f}(\vect{x}^{k}_l,\vect{v}^{k}_l)\right) &&+ \sum_{l=0}^Mqq_{m+1,l}\vect{f}(\vect{x}^{k}_l,\vect{v}^{k}_l)
\end{alignat}
\end{subequations}
for $m=0,\ldots,M-1$ and $k=0,\ldots,K$.
These variables are still vectors, i.\,e.~we have $\vect{v}^{k}_{m},\vect{x}^{k}_{m},\vect{f}(\vect{x}^{k}_l,\vect{v}^{k}_l) \in\mathbb{R}^d$.
Note that the storage required for the $(M+1)\times(M+1)$-matrices~\eqref{eq:Q_matrices} is usually negligible with respect to the size of $d$. 
%For multi-particle setups, $d$ does not only represent the dimension of each particle but also the number of particles in an ensemble.
The following observations are useful to simplify these formulas:
\begin{itemize}
	\item The matrix $\vect{Q}_T$ is a lower diagonal matrix, while $\vect{Q}_x$ is even strictly lower diagonal.
	Thus, the sum over the difference of the function values can be terminated at $m+1$ in~\eqref{eq:node_to_zero_v} and $m$ in~\eqref{eq:node_to_zero_x}. 
	The formula for $\vect{x}_{m+1}^{k+1}$ is therefore fully explicit, the formula $\vect{v}^{k+1}_{m+1}$ is semi-implicit.
    In the following section, this semi-implicit update will be reformulated in an explicit way for problems of the form~\eqref{eq:base_eq} using the Boris ''trick''.
	\item Since $\tau_0 = t_n$, the first row as well as the first column of $\vect{Q}$ is zero. 
	The last sum in both formulas can therefore start at $l=1$, which is also true for the sum over the initial values $\vect{v}_0$ in~\eqref{eq:node_to_zero_x}.
	\item Independently of $q^T_{m+1,0}$ and $q^x_{m+1,0}$, the first term in the summation over the difference of the function values is always zeros, because $\vect{x}_0^k = \vect{x}_0^{k+1} = \vect{x}_0$ and $\vect{v}_0^k = \vect{v}_0^{k+1} = \vect{v}_0$. Thus, these sums can start at $l=1$, too.
	\item The formulations using the initial conditions $\vect{x}_0,\vect{v}_0$ for each node are called ``$0$-to-node'' formulations. 
	More conveniently, this can be reformulated in ``node-to-node'' form using the matrix $\vect{S}$ instead of $\vect{Q}$, where the $m$th row of $\vect{S}=(s_{m,j})_{m,j=0,\ldots,M}$ is defined as the difference between the $m$th and the $(m-1)$th row of $\vect{Q}$ (starting with a row of zeros). 
	The matrices $\vect{S}_x=(s^x_{m,j})_{m,j=0,\ldots,M}$ and $\vect{SQ}=(sq_{m,j})_{m,j=0,\ldots,M}$ are defined analogously to $\vect{Q}_x$ and $\vect{Q}\vect{Q}$.
	Note that the sum over the $(m+1)$th row of $\vect{S}$ is equal to $\Delta\tau_m\coloneq\tau_m-\tau_{m-1}$, which defines the factor in front of $\vect{v}_0$.
\end{itemize}
Based on the comments above, taking the difference between~\eqref{eq:node_to_zero_q} for $m+1$ and $m$ gives
\begin{subequations}\label{eq:sdc_ntn}
    \begin{alignat}{3}
        \vect{v}^{k+1}_{m+1}  &= \vect{v}^{k+1}_m &&+ \frac{\Delta\tau_{m+1}}{2}\left(\vect{f}(\vect{x}^{k+1}_{m+1},\vect{v}^{k+1}_{m+1}) - \vect{f}(\vect{x}^{k}_{m+1},\vect{v}^{k}_{m+1})\right)\notag\\
        & &&+ \frac{\Delta\tau_{m+1}}{2}\left(\vect{f}(\vect{x}^{k+1}_{m},\vect{v}^{k+1}_{m}) - \vect{f}(\vect{x}^{k}_{m},\vect{v}^{k}_{m})\right)
        &+ \sum_{l=1}^Ms_{m+1,l}\vect{f}(\vect{x}^{k}_l,\vect{v}^{k}_l)\label{eq:sdc_ntn_v},\\
        \vect{x}^{k+1}_{m+1} &= \vect{x}^{k+1}_m  + \Delta\tau_{m+1}\vect{v}_0
        &&+ \sum_{l=1}^{m}s^x_{m+1,l}\left(\vect{f}(\vect{x}^{k+1}_l,\vect{v}^{k+1}_l) - \vect{f}(\vect{x}^{k}_l,\vect{v}^{k}_l)\right)
        &+ \sum_{l=1}^Msq_{m+1,l}\vect{f}(\vect{x}^{k}_l,\vect{v}^{k}_l),\label{eq:sdc_ntn_x}
    \end{alignat}
\end{subequations}
for $m=0,\ldots,M-1$ and $k=0,\ldots,K$. This is the ``node-to-node'' formulation of SDC with velocity-Verlet integrator as base method.
Equation~\eqref{eq:sdc_ntn} provides the formulation of Boris-SDC that would actually be implemented: Once values from iteration $k$ are known, the sums involving the quadrature weights can be computed and the step from $m$ to $m+1$ is then essentially a velocity-Verlet step with additional known terms on the right-hand side.

For~\eqref{eq:sdc_ntn}, values for $k=0$ are provided by a simple copy of the initial value to all nodes, see~\eqref{eq:full_matrix_sdc}.
A $k$-times application of these formulas provides approximations $\tvect{v}^k$ and $\tvect{x}^k$ to $\tvect{v}$ and $\tvect{x}$. 
Both are then used to form $\tvect{f}(\tvect{u}^k)$, which in turn serves as input for the collocation formulation to approximate $\vect{u}^k_{n+1}$ via
\begin{align}\label{eq:sdc_coll_end}
\vect{u}^k_{n+1} = \tilde{\vect{C}}_\mathrm{coll}\vect{T}_\mathrm{P}\vect{u}_0 + \tilde{\vect{Q}}_\mathrm{coll}\tvect{f}(\tvect{u}^k),
\end{align}
see Eqs.~\eqref{eq:final_step_coll} and~\eqref{eq:psdc}.
Using $M=2$ Gauss-Lobatto nodes, i.\,e.~$t_n = \tau_0 = \tau_1$ and $\tau_M = \tau_2 = t_{n+1}$, and only a single iteration, the formulas~\eqref{eq:sdc_ntn} yield the standard velocity-Verlet scheme~\eqref{eq:vv}:
A brief calculation shows that for this case 
\begin{align}
    S = \begin{pmatrix}
        0 & 0 & 0\\
        0 & 0 & 0\\
        0 & \frac{\Delta\tau_{2}}{2} & \frac{\Delta\tau_{2}}{2}
    \end{pmatrix},\quad
    S_x = \begin{pmatrix}
    0 & 0 & 0\\
    0 & 0 & 0\\
    0 & \frac{(\Delta\tau_{2})^2}{2} & 0
    \end{pmatrix},\quad
    SQ = \begin{pmatrix}
        0 & 0 & 0\\
        0 & 0 & 0\\
        0 & \frac{(\Delta\tau_{2})^2}{4} & \frac{(\Delta\tau_{2})^2}{4}
        \end{pmatrix}
\end{align}
and $\Delta\tau_1 = 0$.
The collocation formula~\eqref{eq:sdc_coll_end} is obsolete here (since $\tau_M = t_{n+1}$) but nevertheless valid with
\begin{align}
    q = \left(0,\frac{\Delta\tau_2}{2},\frac{\Delta\tau_2}{2}\right),\quad\text{and}\quad Q = S,
\end{align}
see the definitions in~\eqref{eq:tilde_defs}, so that update and iteration formula result in the same expression.
Hence, the first iteration of SDC with velocity-Verlet as base integrator on $M$ Gauss-Lobatto nodes is equivalent to applying velocity-Verlet on each node, as long as the initial value at each node is a copy of the initial value.

\subsection{Boris-SDC}\label{ssec:boris}

For the specific equations of motion under investigation in this work, the right-hand side $\vect{f}$ as stated in~\eqref{eq:rhs_B-field} is given by
\begin{align}
	\vect{f}(\vect{x},\vect{v}) = \alpha\left[\vect{E}(\vect{x},t) + \vect{v}\times\vect{B}\right]\label{eq:rhs_boris}
\end{align}
with a constant magnetic field $\vect{B}(\vect{x},t) = \vect{B}$ for simplicity (see note on non-constant $\vect{B}$-fields at the end of this section). 
Thus, the update formulas for the velocity component $\vect{v}$ in the standard velocity-Verlet integrator~\eqref{eq:vv_v} as well as in the SDC iteration~\eqref{eq:sdc_ntn_v} seem to require the solution of an implicit system.
For velocity-Verlet integration, the state-of-the-art approach is the Boris integration method~\cite{Boris1970,Birdsall1985}.
Here, \eqref{eq:vv_v} is rewritten as
\begin{align}
	\frac{\vect{v}_{n+1}-\vect{v}_n}{\Delta t} = \alpha\left[\vect{E}_{n+\frac{1}{2}} + \frac{\vect{v}_{n+1}+\vect{v}_n}{2}\times\vect{B}\right].\label{eq:rhs_B-field_findiff}
\end{align}
The half-step subscript corresponds to the average electric field at times $t_n$ and $t_{n+1}$, i.\,e.~$\vect{E}_{n+\frac{1}{2}} \coloneq \frac{1}{2}\left(\vect{E}(\vect{x}_n,t_n) + \vect{E}(\vect{x}_{n+1},t_{n+1})\right)$.
The idea of the Boris integrator is to separate the electric and magnetic forces.
To this end, we define
\begin{align}
	\vect{v}^- \coloneq \vect{v}_{n} + \frac{\alpha\Delta t}{2}\vect{E}_{n+\frac{1}{2}}
  \quad\text{and}\quad
	\vect{v}^+ \coloneq \vect{v}_{n+1} - \frac{\alpha\Delta t}{2}\vect{E}_{n+\frac{1}{2}}, \label{eq:Borisvpandm}
\end{align}
so that
\begin{align}
	\frac{\vect{v}^+-\vect{v}^-}{\Delta t} &= \frac{\alpha}{2}(\vect{v}^++\vect{v}^-)\times\vect{B}, \label{eq:Borisvpm}
\end{align}
which can be shown to correspond to a simple rotation, i.\,e.~$|\vect{v}^+| = |\vect{v}^-|$ and can be solved for $\vect{v}^+$ explicitly (cf.~\cite{Birdsall1985}) using $\vect{v}^+ = \vect{v}^- + \left(\vect{v}^-+\vect{v^-}\times\vect{t}\right)\times\vect{s}$ with $\vect{t} = \alpha\vect{B}\cdot{\Delta t}/{2}$ and $\vect{s} = 2\vect{t}/\left(1+|\vect{t}|^2\right)$. % see \cite[Chapter~15--5]{Birdsall1985} for details on time-dependent magnetic fields. Essentially, they are no big deal.
Using~\eqref{eq:Borisvpandm}, the new velocity $\vect{v}_{n+1}$ can therefore be computed explicitly.
A relativistic generalization of this method is straightforward~\cite[Chapter~15--4]{Birdsall1985}.

We can use and extend this idea to resolve the implicit dependency in the SDC iteration~\eqref{eq:sdc_ntn_v}.
More precisely, we rewrite the (seemingly implicit) update for the $(m+1)$th component of $\tvect{v}$ as
\begin{align}
\begin{split}
\frac{\vect{v}^{k+1}_{m+1}  - \vect{v}^{k+1}_m}{\Delta\tau_{m+1}} = &\frac{1}{2}\left(\vect{f}(\vect{x}^{k+1}_{m+1},\vect{v}^{k+1}_{m+1}) + \vect{f}(\vect{x}^{k+1}_{m},\vect{v}^{k+1}_{m})\right)
-\frac{1}{2}\left(\vect{f}(\vect{x}^{k}_{m+1},\vect{v}^{k}_{m+1}) + \vect{f}(\vect{x}^{k}_{m},\vect{v}^{k}_{m})\right)\\
&+\frac{1}{\Delta\tau_{m+1}}\sum_{l=1}^Ms_{m+1,l}\vect{f}(\vect{x}^{k}_l,\vect{v}^{k}_l).
\end{split}
\end{align}
Note that the second and third summand of the right-hand side of this equation only depend on values at iteration $k$, i.\,e.~these summands have been computed in the previous iteration. 
We define
\begin{align}
	\Delta\tau_{m+1}\vect{c}^k \coloneq -\frac{\Delta\tau_{m+1}}{2}\left(\vect{f}(\vect{x}^{k}_{m+1},\vect{v}^{k}_{m+1}) + \vect{f}(\vect{x}^{k}_{m},\vect{v}^{k}_{m})\right) +\sum_{l=1}^Ms_{m+1,l}\vect{f}(\vect{x}^{k}_l,\vect{v}^{k}_l),
\end{align}
so that
\begin{align}
	\frac{\vect{v}^{k+1}_{m+1}  - \vect{v}^{k+1}_m}{\Delta\tau_{m+1}} = \frac{\vect{f}(\vect{x}^{k+1}_{m+1},\vect{v}^{k+1}_{m+1})
		+ \vect{f}(\vect{x}^{k+1}_{m},\vect{v}^{k+1}_{m})}{2} + \vect{c}^k.
\end{align}
With the particular right hand side~\eqref{eq:rhs_boris}, this yields
\begin{align}
	\frac{\vect{v}_{m+1}^{k+1}-\vect{v}_m^{k+1}}{\Delta\tau_{m+1}} = \alpha\left[\vect{E}_{m+\frac{1}{2}}^{k+1} + \frac{\vect{v}_{m+1}^{k+1}+\vect{v}_m^{k+1}}{2}\times\vect{B}\right] + \vect{c}^k,
\end{align}
Except for the $\vect{c}^k$-term (which is known from the previous iteration) this has the very same structure as Equation~\eqref{eq:rhs_B-field_findiff}.
Extending the idea of~\eqref{eq:Borisvpandm}, we define
\begin{align}
\vect{v}^- \coloneq \vect{v}_m^{k+1} + \frac{\Delta\tau_{m+1}}{2}\left(\alpha\vect{E}^{k+1}_{m+\frac{1}{2}} + \vect{c}^k\right)
\quad\text{and}\quad
\vect{v}^+ \coloneq \vect{v}_{m+1}^{k+1} - \frac{\Delta\tau_{m+1}}{2}\left(\alpha\vect{E}^{k+1}_{m+\frac{1}{2}} + \vect{c}^k\right)\label{eq:Borisvp_sdc}
\end{align}
to obtain
\begin{align}
\frac{\vect{v}^+-\vect{v}^-}{\Delta\tau_{m+1}} = \frac{\alpha}{2}(\vect{v}^++\vect{v}^-)\times\vect{B},
\end{align}
which is precisely of the type of~\eqref{eq:Borisvpm}. 
As noted before, this can be solved explicitly for $\vect{v}^+$, so that~\eqref{eq:Borisvp_sdc} can be used to determine $\vect{v}_{m+1}^{k+1}$.
This gives us an explicit solver for the seemingly implicit SDC update~\eqref{eq:sdc_ntn_v} and can be implemented directly into an existing SDC algorithm without further modifications.

We note that this approach can be easily extended to non-constant magnetic fields $\vect{B}(\vect{x},t)$ as follows. 
Instead of Eq.~\eqref{eq:rhs_B-field_findiff} we now have
\begin{align}
	\frac{\vect{v}_{n+1}-\vect{v}_n}{\Delta t} &= \alpha\left[\vect{E}_{n+\frac{1}{2}} + \frac{1}{2}\vect{v}_{n+1}\times\vect{B}_{n+1} + \frac{1}{2}\vect{v}_{n}\times\vect{B}_{n}\right]\\
	&= \alpha\left[\vect{E}_{n+\frac{1}{2}} + \frac{\vect{v}_{n+1}+\vect{v}_{n}}{2}\times\vect{B}_{n+1}\right] + \frac{\alpha}{2}\vect{v}_n\times(\vect{B}_n-\vect{B}_{n+1}).
\end{align}
The first part has again the same structure as Eq.~\eqref{eq:rhs_B-field_findiff}, while the last part does not depend on $\vect{v}_{n+1}$ and can thus be treated separately, i.\,e. as part of the $\vect{c}^k$-term in~\eqref{eq:Borisvp_sdc}.

\section{Numerical results}\label{sec:numerics}
To evaluate the numerical properties of the Boris-SDC integrator, we study particles in a standard Penning trap~\cite{Penning1936}.
Being confined to a limited volume due to an external magnetic and electric field, the particles' characteristic properties, such as trajectories in real and phase space, energy conservation, and stability of the integration scheme, can conveniently be analyzed.
We consider both the case of a single particle, where an analytic reference solution for the particle's trajectory is available, as well as the case of many particles.
All Boris-SDC runs use Gauss-Lobatto quadrature nodes, see the discussion in~\ref{sec:collprop}.

We follow the analysis in~\cite{Patacchini2009} and choose a constant magnetic field $\vect{B}=\frac{\omega_B}{\alpha}\cdot\hat{\vect{e}}_z$ along the $z$-axis with the particle's charge-to-mass ratio $\alpha=\nicefrac{q}{m}$ so that
\begin{align}
   \vect{v}\times\vect{B} = \frac{\omega_B}{\alpha} 
   \begin{pmatrix}
      0 & 1 & 0 \\
      -1 & 0 & 0\\
      0 & 0 & 0
   \end{pmatrix}\vect{v} .\label{eq:Bpenning}
\end{align} 
The electric field $\vect{E}(\vect{x}_i) = \vect{E}_\text{ext}(\vect{x}_i) + \vect{E}_\text{int}(\vect{x}_i)$ experienced by a particle at position $\vect{x}_i$ is composed of an ideal quadrupole potential distribution leading to 
\begin{align}
   \vect{E}_\text{ext}(\vect{x}_i) = -\epsilon\frac{\omega_E^2}{\alpha}
   \begin{pmatrix}
      1 & 0 & 0 \\
      0 & 1 & 0 \\
      0 & 0 & -2
   \end{pmatrix}\vect{x}_i\label{eq:Epenning}
\end{align}
and the inter-particle Coulomb interaction (cgs units)
\begin{align}
  \vect{E}_\text{int}(\vect{x}_i) = \sum_{\begin{subarray}{c}k=1\\k\neq i\end{subarray}}^{N_\text{particles}}\mkern-8mu Q_k \frac{\vect{x}_i-\vect{x}_k}{(|\vect{x}_i-\vect{x}_k|^2+\lambda^2)^{\nicefrac{3}{2}}}.
\end{align}
To avoid numerical heating due to spurious close encounters, we follow the standard approach and regularize the Coulomb pole here by means of the smoothing parameter $\lambda>0$.

In the case $\epsilon=-1$, particles are confined in $z$-direction due to the attractive nature of the potential.
This setup corresponds to a Penning trap configuration.
For $\epsilon=+1$, the particle will escape along the $z$-axis but -- due to the magnetic field -- will still be confined to a limited orbit in $x$- and $y$-direction.

%%%%%%%%%%%%%%%%%%%%%%%%%%%%%%%%%%%%%%%%%%%%%%%%%%%%%%%%%%%%%%%%%%%%%%%%%%%%%%%%
%%%%%%%%%%%%%%%%%%%%%%%%%%%%%%%%%%%%%%%%%%%%%%%%%%%%%%%%%%%%%%%%%%%%%%%%%%%%%%%%
\subsection{Single particle in a Penning trap}

\begin{figure}[tb]
    \begin{minipage}{0.4\textwidth}
    \centering
  \begin{tabular}{lc}
    \toprule
      $\alpha$        & $1.0$ \\
      $t_\text{end}$  & $16.0$ \\
      $\vec{x}(t=0)$  & $( 10, 0,   0)^T$ \\
      $\vec{v}(t=0)$  & $(100, 0, 100)^T$ \\
      $\omega_E$      & $4.9$  \\
      $\omega_B$      & $25.0$ \\
      $N_\text{steps}$ & variable \\
    \bottomrule
  \end{tabular}
  \captionof{table}{Setup parameters chosen for the case of a single classical particle in the Penning trap.}
  \label{tab:param1}
  \end{minipage}\hfill
  \begin{minipage}{0.49\textwidth}
  \centering
  \includegraphics[width=\textwidth]{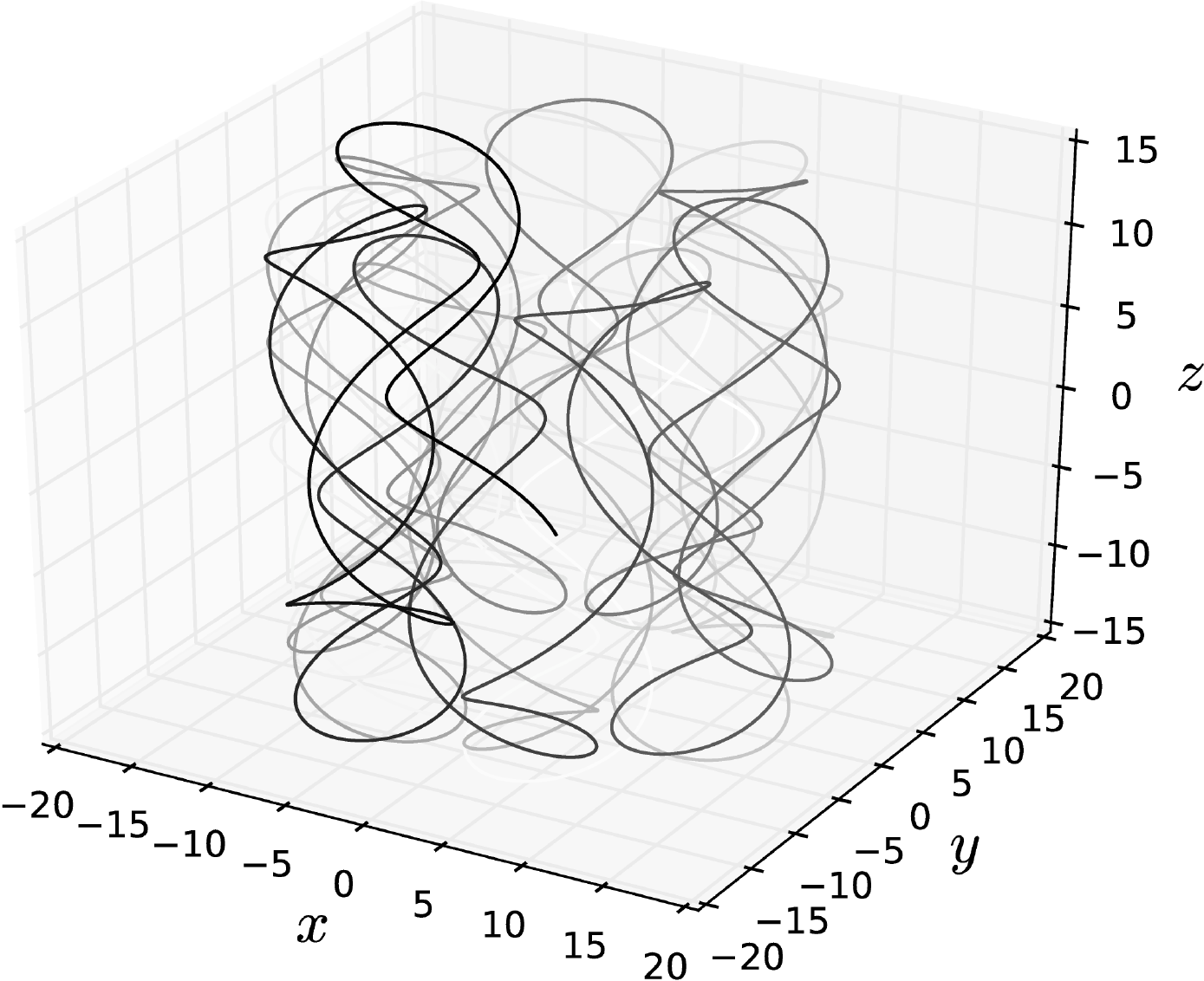}
  \captionof{figure}{Particle trajectory for the parameters of Table~\ref{tab:param1}.
  Evolution in time is indicated by the line's changing shade from light ($t=0$) to dark ($t=t_\text{end}$).}
  \label{fig:trajectory}
  \end{minipage}
\end{figure}

Solving the equations of motion~\eqref{eq:newton} with the magnetic field~\eqref{eq:Bpenning} and the electric field~\eqref{eq:Epenning}
for a single particle inside the Penning trap is a standard textbook task, see e.\,g.~\cite{Brown1986,Werth2009}.
Movement in $z$ direction is a harmonic oscillation, decoupled from the other coordinates,
\begin{align}
  %\begin{split}
    z(t) &= z(0)\cos(\tilde{\omega}t) + \frac{v_z(0)}{\tilde{\omega}}\sin(\tilde{\omega}t), & \tilde{\omega}&\coloneq\sqrt{-2\epsilon}\cdot\omega_E. \label{eq:omegatilde}
    %\frac{d}{\d t}z(t) &= -\sqrt{2}\omega_E z_0 \sin(\sqrt{2}\omega_E t) + v_{z,0}\cos(\sqrt{2}\omega_E t) \\
    %\\
    %z_0 &= z(t=0),  \quad  v_{z,0} = v_z(t=0).
  %\end{split}
\end{align}
Here, the role of $\epsilon$ stated before becomes obvious:
For $\epsilon=+1$, the frequency $\tilde{\omega}$ is purely imaginary and the $z$-trajectory diverges, for $\epsilon=-1$ the frequency $\tilde{\omega}\in\mathbb{R}$ and the trajectory $z(t)$ corresponds to a harmonic oscillation.
With the definition $w(t)\coloneq  x(t) + \ii y(t)$, the particle movement in the $x$-$y$ plane is given by
\begin{align}
  w(t) &= (\scrR_+ + \ii\scrI_+)\exp{-\ii\Omega_+t} 
        + (\scrR_- + \ii\scrI_-)\exp{-\ii\Omega_-t}
        &
  \Omega_\pm &\coloneq \frac{1}{2}\left(\omega_B \pm \sqrt{\omega_B^2 + 4\epsilon\omega_E^2}\right)
\end{align}
\begin{align}
  \scrR_- &\coloneq \frac{\Omega_+x(0)+v_y(0)}{\Omega_+-\Omega_-},\quad 
  &
  \scrR_+ &\coloneq x(0) - \scrR_-,
  &
  \scrI_- &\coloneq \frac{\Omega_+y(0)-v_x(0)}{\Omega_+-\Omega_-},
  &
  \scrI_+ &\coloneq y(0) - \scrI_- .
\end{align}
Note that for $\omega_B^2 < -4\epsilon\omega_E^2$ (which can only happen for $\epsilon=-1$), we have $\Omega_\pm\notin\mathbb{R}$ for the revolution frequency.
In this case, the physical setup is unstable and the particle escapes from the trap due to a too weak magnetic or too strong electric field.
Table~\ref{tab:param1} lists the physical parameters used in this section; Figure~\ref{fig:trajectory} shows a visualization of the particle's analytical trajectory.

%%%%%%%%%%%%%%%%%%%%%%%%%%%%%%%%%%%%%%%%%%%%%%%%%%%%%%%%%%%%%%%%%%%%%%%%%%%%%%%%
\subsubsection{Stability}

\begin{figure*}[p]
    \centering
    \begin{minipage}{\figurewidthWW\textwidth}
      \centerline{Classical Boris}
      \includegraphics[width=\textwidth]{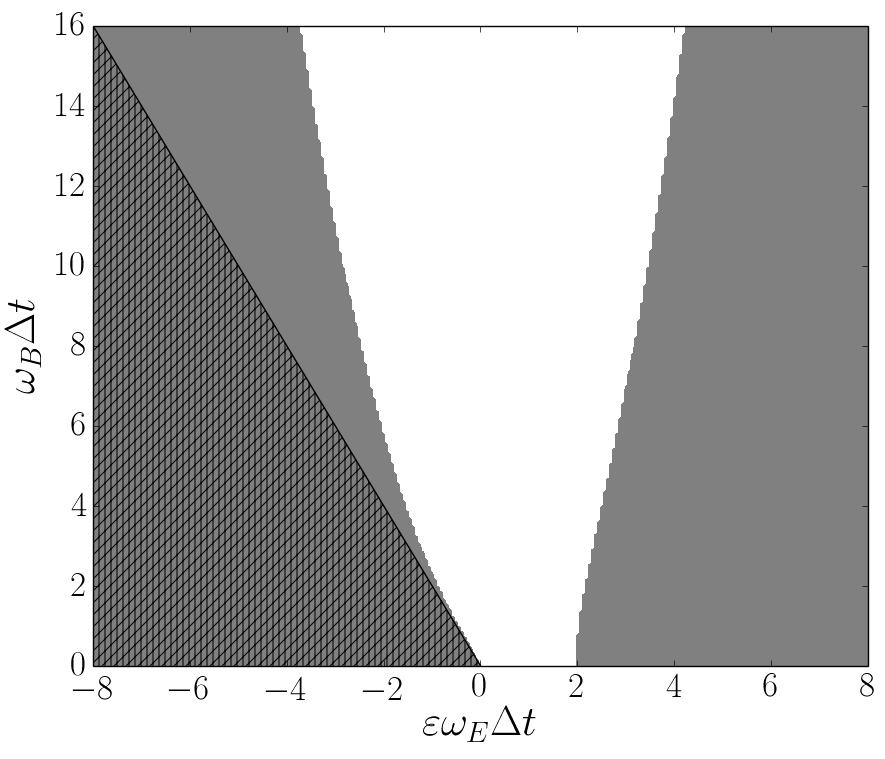}
    \end{minipage}
    \begin{minipage}{\figurewidthWW\textwidth}
        \centerline{Collocation, $M=3$}
        \includegraphics[width=\textwidth]{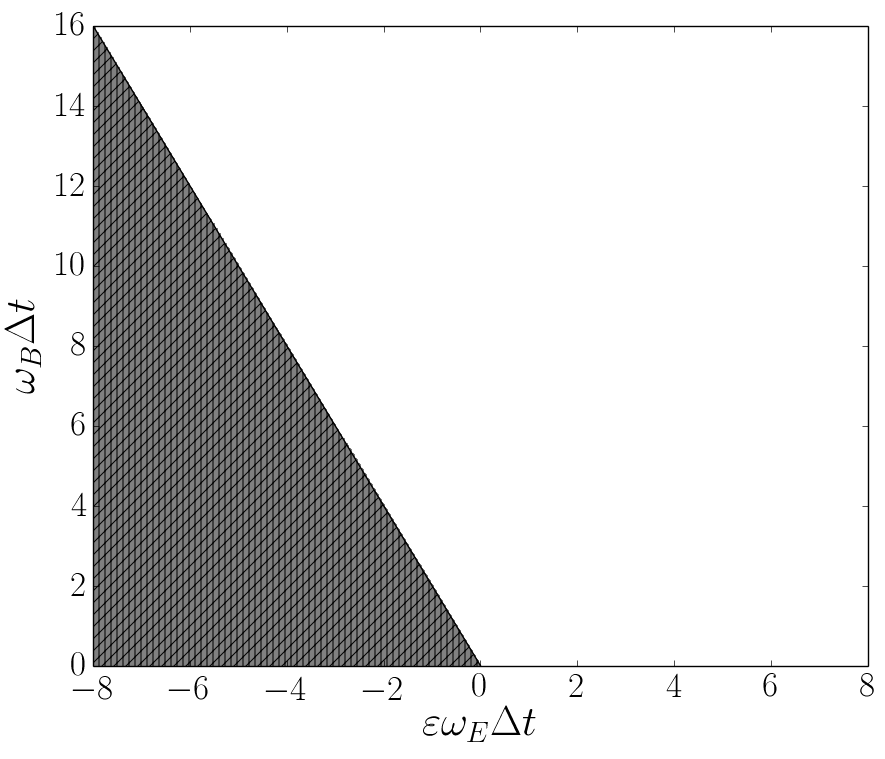}
    \end{minipage}
    \caption{Stability regions for the single-particle system in an ideal Penning trap.
    The dark gray, hatched regions denote the trap's physical instability region $\omega_B^2 < -4\epsilon\omega_E^2$ while light gray areas indicate numerical instability. For $M=5$ nodes (not shown), the plot for the collocation method is identical to the one with $M=3$.}
    \label{fig:stability}
\end{figure*}
\begin{figure*}[p]
      \centering
      \begin{minipage}{\figurewidthWW\textwidth}
          \centerline{$M=3$}
          \includegraphics[width=\textwidth]{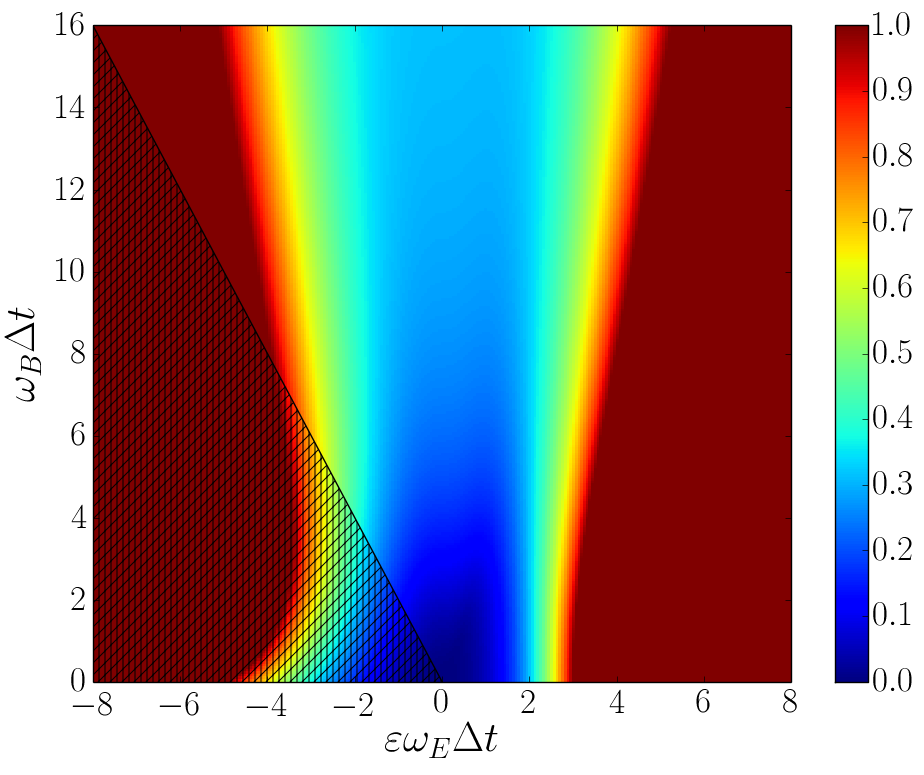}
      \end{minipage}
      \begin{minipage}{\figurewidthWW\textwidth}
          \centerline{$M=5$}
          \includegraphics[width=\textwidth]{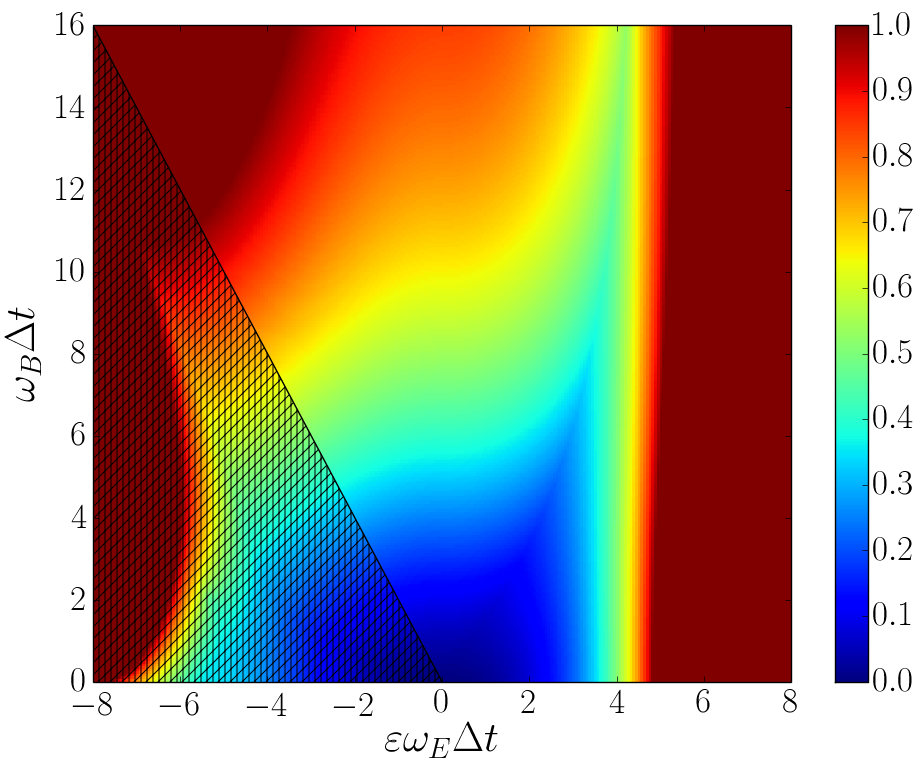}
      \end{minipage}
      \caption{Convergence regions for the single-particle system in an ideal Penning trap.
          Colors encode the spectral radius of the SDC iteration matrix $\vect{K}_\mathrm{sdc}$ with $M=3, 5$ nodes and sufficiently many iterations to reach a residual tolerance of $10^{-12}$.
          The dark gray, hatched region denotes the trap's physical instability region $\omega_B^2 < -4\epsilon\omega_E^2$.}
      \label{fig:convergence}
\end{figure*}
\begin{figure*}[p]
      \centering
        \begin{minipage}{\figurewidthWW\textwidth}
          \centerline{$M=3$}
          \includegraphics[width=\textwidth]{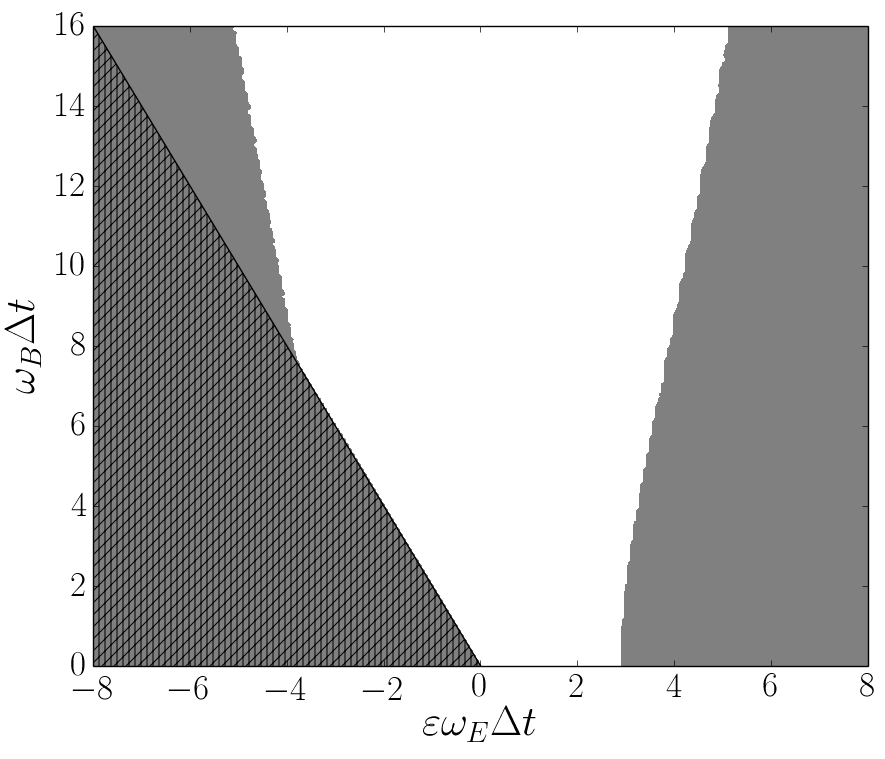}
      \end{minipage}
      \begin{minipage}{\figurewidthWW\textwidth}
          \centerline{$M=5$}
          \includegraphics[width=\textwidth]{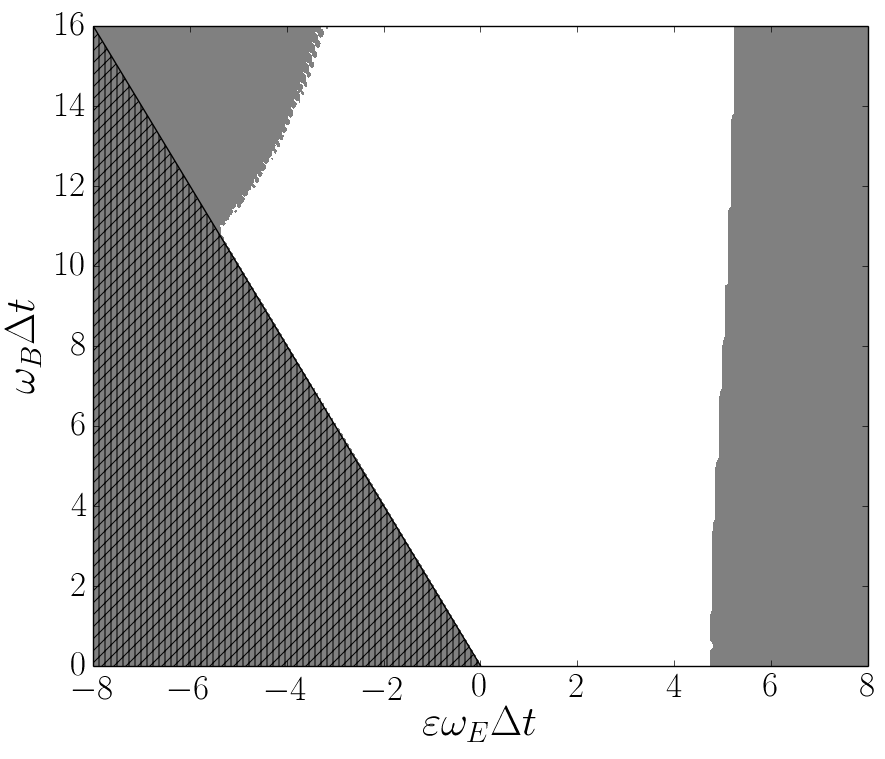}
      \end{minipage}
      \caption{Stability regions for Boris-SDC for the single-particle system in an ideal Penning trap.
      	Sufficiently many iterations are performed to reach a residual tolerance of $10^{-12}$.
          Light gray color indicates numerical instability of the method using $\vect{P}_\mathrm{sdc}$ with $M=3,5$ nodes.
          The dark gray, hatched region denotes the trap's physical instability region $\omega_B^2 < -4\epsilon\omega_E^2$.}
      \label{fig:sdcstability}
  \end{figure*}  

As a first study, we analyze the stability of the classical Boris integrator, the collocation method, and Boris-SDC for the Penning trap and discuss the relation of the stability of Boris-SDC to the convergence of the corresponding SDC iteration.
Note that the numerical stability of SDC has been extensively studied for the generic test equation with a single complex eigenvalue in~\cite{DuttEtAl2000}.
For a modified version of the original Boris scheme, a numerical stability analysis for an ideal Penning trap has been performed in~\cite{Patacchini2009}.

We assess stability by the largest absolute value of the eigenvalues, that is the spectral radius, of the method's update matrix. For Boris-SDC this is the matrix $\tilde{\vect{P}}^k_\mathrm{sdc}$ defined in~\eqref{eq:psdc} and for the collocation method $\tilde{\vect{P}}_\mathrm{coll}$ defined in~\eqref{eq:coll_update}.
A method is stable for a specific configuration, if and only if the spectral radius is smaller than or equal to unity.

Figure~\ref{fig:stability} shows the resulting stability region for the classical Boris integrator and a collocation method with $M=3$.
The hatched, dark gray area on the left indicates the physical instability of the system, where $\omega_B^{2} < -4\epsilon\omega_E^{2}$ so that the magnetic field is too weak to confine the particle.
For the classical Boris integrator, there is in addition a zone of numerical instability (light gray), where the method is not stable although the physical problem already is.
Furthermore, a significant area of numerical instability is present for values of $\varepsilon \omega_E \Delta t > 2$ on the right side.
In contrast, the stability domain of the collocation method is identical to the domain where the physical setup is stable:
There is no region of additional numerical instability.
For $M=5$, the stability domain of the collocation method is identical (not shown).
The collocation method is hence stable for every physically stable configuration of the single-particle Penning trap.   

\paragraph{Convergence of Boris-SDC iteration} For Boris-SDC, the convergence properties of the iteration computing the collocation solution have significant impact on the stability:
In cases where the iteration converges poorly or not all, we cannot expect to recover the stability properties of the collocation solution.
Convergence of the SDC iteration is governed by the spectral radius of the SDC iteration matrix $\vect{K}_\mathrm{sdc}$ (not to be confused with the Boris-SDC update matrix).
If and only if the spectral radius $\spectral(\vect{K}_\mathrm{sdc}) < 1$, the iteration will ultimately converge and the norm of the residual~\eqref{eq:residual} will go to zero.
For values close to unity, however, convergence can be unfeasibly slow, resulting in a very large number of required iterations.
Figure~\ref{fig:convergence} shows the spectral radius of $\vect{K}_\mathrm{sdc}$ for $M=3$ and $M=5$ collocation nodes.
Small values (blue) indicate fast convergence, values close to unity (yellow and light red) slow convergence and values larger than unity (dark red) indicate divergence of the SDC iteration.
For $M=3$, the area where Boris-SDC shows good convergence does roughly coincide with the stability domain of the classical Boris integrator.
Interestingly, for small values of $\omega_B  \Delta t$, Boris-SDC also converges well in the region of physical instability.
For $M=5$ nodes, the picture is quite different. 
Here, Boris-SDC has a somewhat larger convergence domain than with three nodes for small values of $\omega_B \Delta t$.
As $\omega_B \Delta t $ increases, however, independently of the value for $\omega_E  \Delta t$, convergence eventually starts to deteriorate.
This suggests that if the particle motion induced by the magnetic field becomes strongly under-resolved, Boris-SDC fails to converge for larger values of $M$.
Strategies exist to improve convergence of SDC, see e.\,g.~\cite{HuangEtAl2006,Weiser2014}, but studying their effect on the convergence of Boris-SDC is left for future work.
Preliminary tests suggest that these strategies can be applied here as well.

\paragraph{Stability of Boris-SDC} Figure~\ref{fig:sdcstability} shows the stability domain of Boris-SDC, that is the region where the spectral radius of the update matrix $\tilde{\vect{P}}^{k}_\mathrm{sdc}$ is less  or equal than one.
The iteration count $k$ is chosen to satisfy a residual tolerance of $r\leq10^{-12}$.
Thus, in the regions where Boris-SDC is converging slowly according to Figure~\ref{fig:convergence}, more iterations are performed.
Because the underlying collocation method is stable outside the region of physical instability, numerical instabilities of Boris-SDC can only arise due to failure of the iteration to converge.
The regions of numerical instability (light gray) for Boris-SDC at the right and on the upper left in Figure~\ref{fig:sdcstability} correspond to the dark red regions in Figure~\ref{fig:convergence} where $\spectral(\vect{K}_\mathrm{sdc})>1$.
For a spectral radius smaller than one, the iteration will converge, although probably slowly, and eventually recover the stability of the underlying collocation scheme. 
Although Boris-SDC does not fully maintain the stability properties of the underlying collocation method, the stability regions for both $M=3$ and $M=5$ are larger than for the classical Boris-integrator.

%%%%%%%%%%%%%%%%%%%%%%%%%%%%%%%%%%%%%%%%%%%%%%%%%%%%%%%%%%%%%%%%%%%%%%%%%%%%%%%%
\subsubsection{Order of convergence}
A key advantage of SDC-based integration methods is their ability to easily generate methods of high order.
Typically, for an Euler base method, each iteration or sweep increases the order by one, up to the order of the underlying collocation method, see e.\,g.~\cite{ShuEtAl2007}.
If an order-$p$ base method is used in combination with equidistant quadrature nodes, each iteration increases the overall order by~$p$~\cite{ChristliebEtAl2010}.
This property, however, does {\it not} necessarily hold if non-equidistant nodes like e.\,g.~Gauss nodes are used.
In the examples presented here, however, we generally observe an increase of the overall order by two with each sweep of the second order Boris integrator.

  \begin{figure}[t]
      \centering
      \begin{minipage}{\figurewidthTR\textwidth}
          \centerline{$M=3$}
          \includegraphics[width=\textwidth]{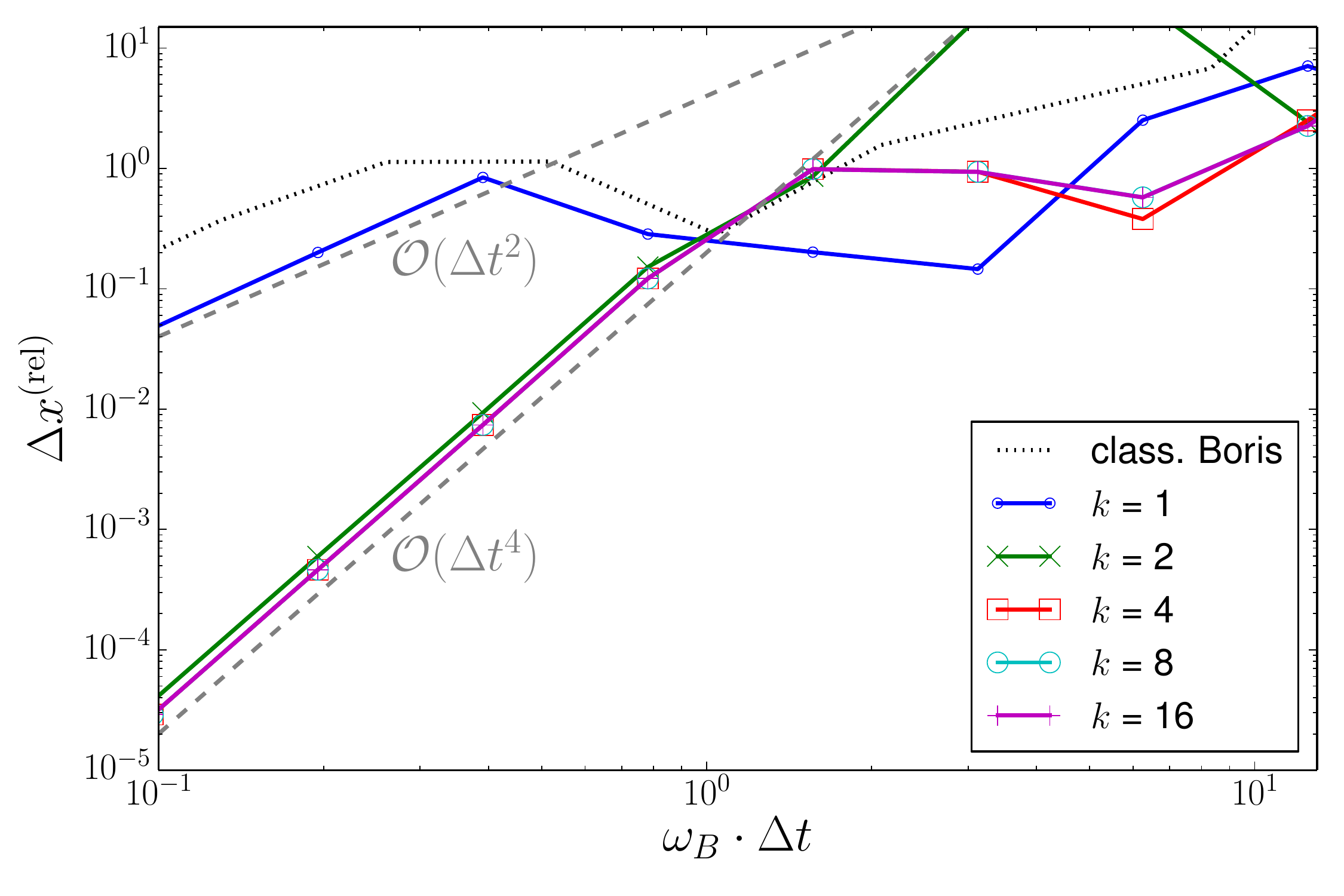}
      \end{minipage}
      \begin{minipage}{\figurewidthTR\textwidth}
          \centerline{$M=5$}
          \includegraphics[width=\textwidth]{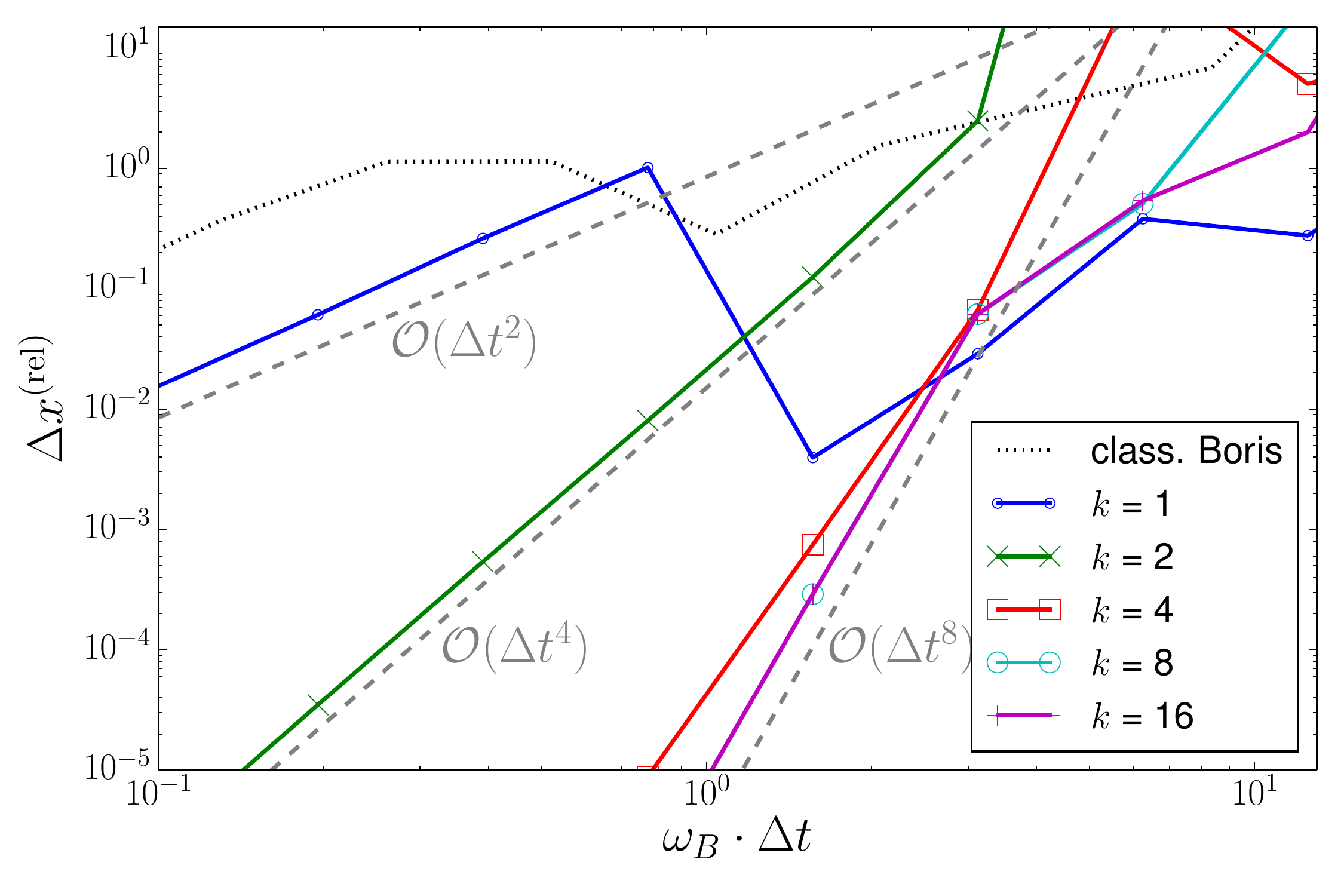}
      \end{minipage}
    \caption{Relative error $\Delta x^{(\text{rel})}$, see~\eqref{eq:rel_error_x}, for the $x$ coordinate of the particle's final position in the Penning trap as a function of time step size for 3 and 5 Gauss-Lobatto collocation nodes per time step and different fixed numbers of iterations per SDC sweep.
    As a guide to the eye, dashed lines indicate second, fourth and eighth order and for comparison, the classical Boris method's convergence is shown.
}
    \label{fig:order_GLob}
  \end{figure}

Figure~\ref{fig:order_GLob} shows the relative error 
\begin{align}
	\label{eq:rel_error_x}
	\Delta x^{(\text{rel})}\coloneq \frac{ \left| x-x^{(\text{analyt})} \right| }{ \left| x^{(\mathrm{analyt})} \right| } 
\end{align}
at $t=t_\mathrm{end}$ in the $x$ coordinate of the Boris-SDC method against the analytical solution depending on the length of the time step for $M=3$ and $M=5$ Gauss-Lobatto nodes and different numbers of iterations.
Below, we always report the error in the $x$ coordinate of the position, but the $y$ and $z$ components as well as the velocities show analogous behavior.

The order of the underlying quadrature rule is $2M - 2$, i.\,e.~four for $M=3$ and eight for $M=5$.
As a guide to the eye, lines indicating second, fourth and eighth order are included.
For $M=3$, a single sweep already yields a second-order method.
Two or more sweeps are sufficient to reproduce the convergence order of the underlying fourth order collocation method:
The lines for two and more sweeps are essentially identical.
Also for $M=5$, each iteration raises the order by two.
With four iterations, the eighth order of the underlying collocation method is reached, but going to eight sweeps still yields a small improvement.
Note that the classical Boris integrator corresponds to Boris-SDC with a single iteration (provided the initial values are set adequately as discussed in~\ref{sec:vv_sdc}).
Since the step size $\dt$ is the same for Boris-SDC and the classical integrator, the use of three and five nodes leads to better accuracies for Boris-SDC with $k=1$, though.
%Note that, as discussed in~\ref{sec:energy}, retrieving the symplecticness of the collocation method requires more sweeps than recovering only its order of convergence.

%%%%%%%%%%%%%%%%%%%%%%%%%%%%%%%%%%%%%%%%%%%%%%%%%%%%%%%%%%%%%%%%%%%%%%%%%%%%%%%%
\subsubsection{Residual Control}\label{sec:residual_control}
Considering SDC as an iterative solver for the collocation problems allows to prescribe some tolerance and to iterate until the norm of the residual $r$ defined in~\eqref{eq:residual} is below this tolerance instead of prescribing a fixed number of iterations.
This allows for easy tuning of  precision against performance instead of just fixing a specific convergence order.

\begin{figure}[t]
      \centering
      \begin{minipage}{\figurewidthTR\textwidth}
          \centerline{$M=3$}
          \includegraphics[width=\textwidth]{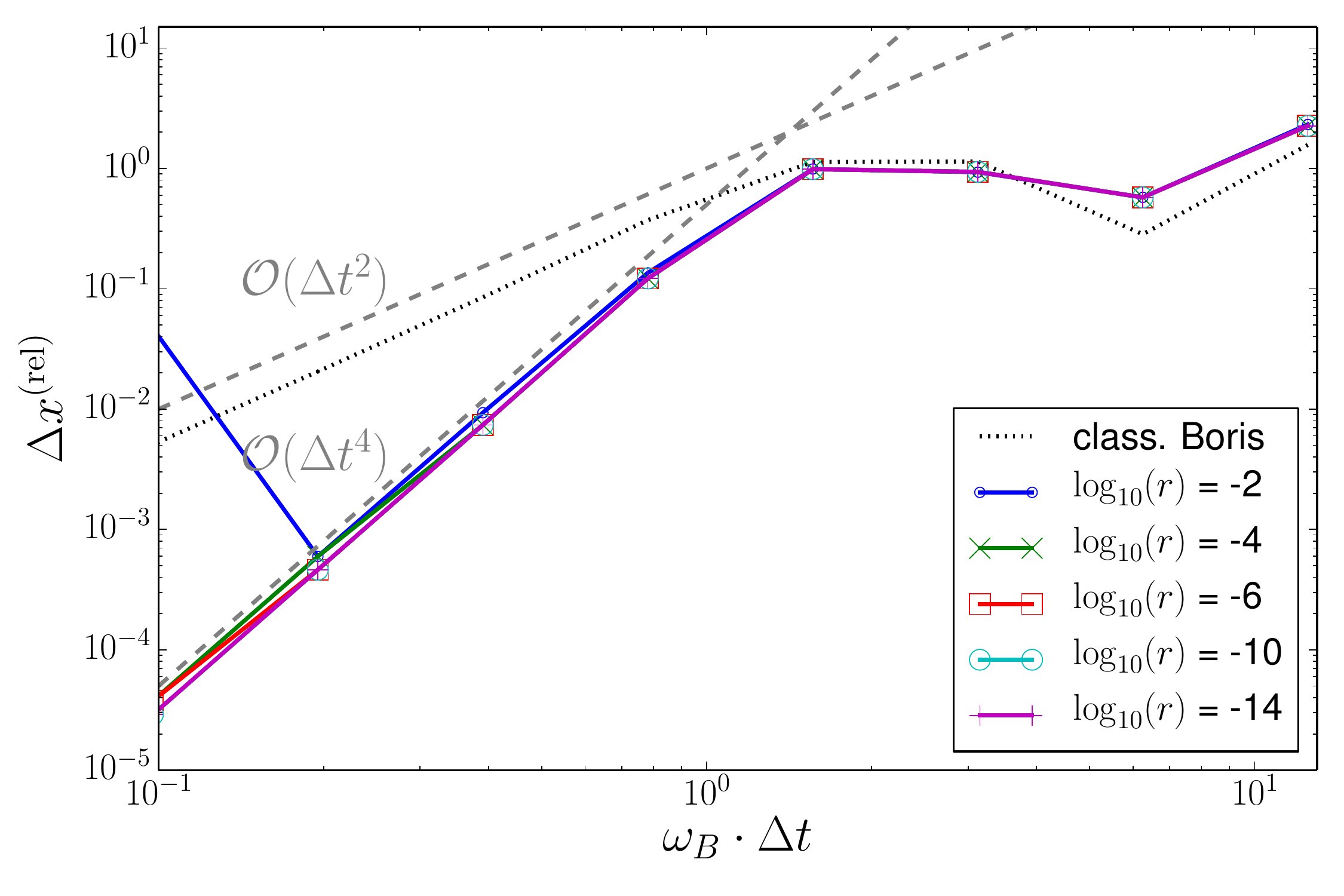}
      \end{minipage}
      \begin{minipage}{\figurewidthTR\textwidth}
          \centerline{$M=5$}
          \includegraphics[width=\textwidth]{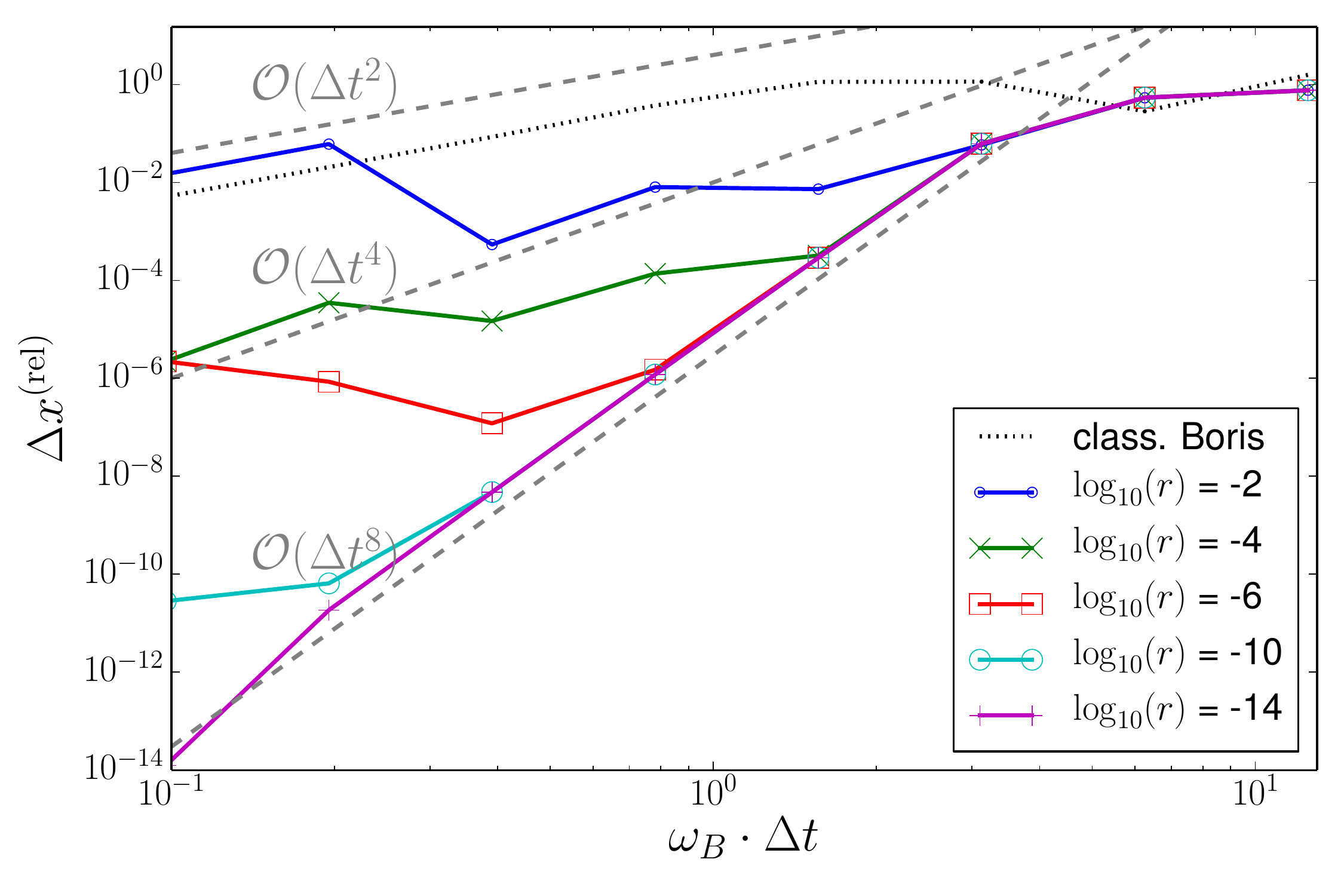}
      \end{minipage}
  \caption{Relative error $\Delta x^{(\text{rel})}$, see~\eqref{eq:rel_error_x}, for the $x$ coordinate of the particle's final position in the Penning trap as a function of time step for 3 and 5 collocation nodes per time step and different tolerances for the residual of the SDC iteration.}
  \label{fig:residual}
\end{figure}
Figure~\ref{fig:residual} shows the relative error $\Delta x^{(\text{rel})}$ of Boris-SDC depending on the time step, again for $M=3$ and $M=5$ Gauss-Lobatto nodes.
Instead of a fixed number of iterations as for Figure~\ref{fig:order_GLob}, the runs here use a prescribed tolerance for the residual of SDC:
On each time step, sweeps are performed until the requested tolerance is met, see also~\cite{SpeckEtAl2014_BIT}.
As the time step is decreased, $\Delta x^{(\text{rel})}$ first decreases with the order of the underlying collocation method unless saturating at a level determined by the prescribed tolerance.
The smaller the tolerance is set, the later the error is saturating.
For a sufficiently large number of quadrature nodes, the error can essentially be brought down to almost machine precision, if desired.

Note that although the residual condition is a per time step measure, in the case studied here the tolerance set for the residual gives a decent indication of the actual final error. 
The values where the error saturates in all experiments closely match the set residual tolerance:
An $r\leq10^{-2}$ tolerance for the residual also results in approximately a $\Delta x^{(\mathrm{rel})}\approx10^{-2}$ final error for example, with the notable exception of a $10^{-4}$ tolerance which results in a significantly smaller error.

%%%%%%%%%%%%%%%%%%%%%%%%%%%%%%%%%%%%%%%%%%%%%%%%%%%%%%%%%%%%%%%%%%%%%%%%%%%%%%%%
\subsubsection{Work--precision}
\begin{figure}[t]
      \centering
      \begin{minipage}{\figurewidthTR\textwidth}
          \centerline{$M=3$}
          \includegraphics[width=\textwidth]{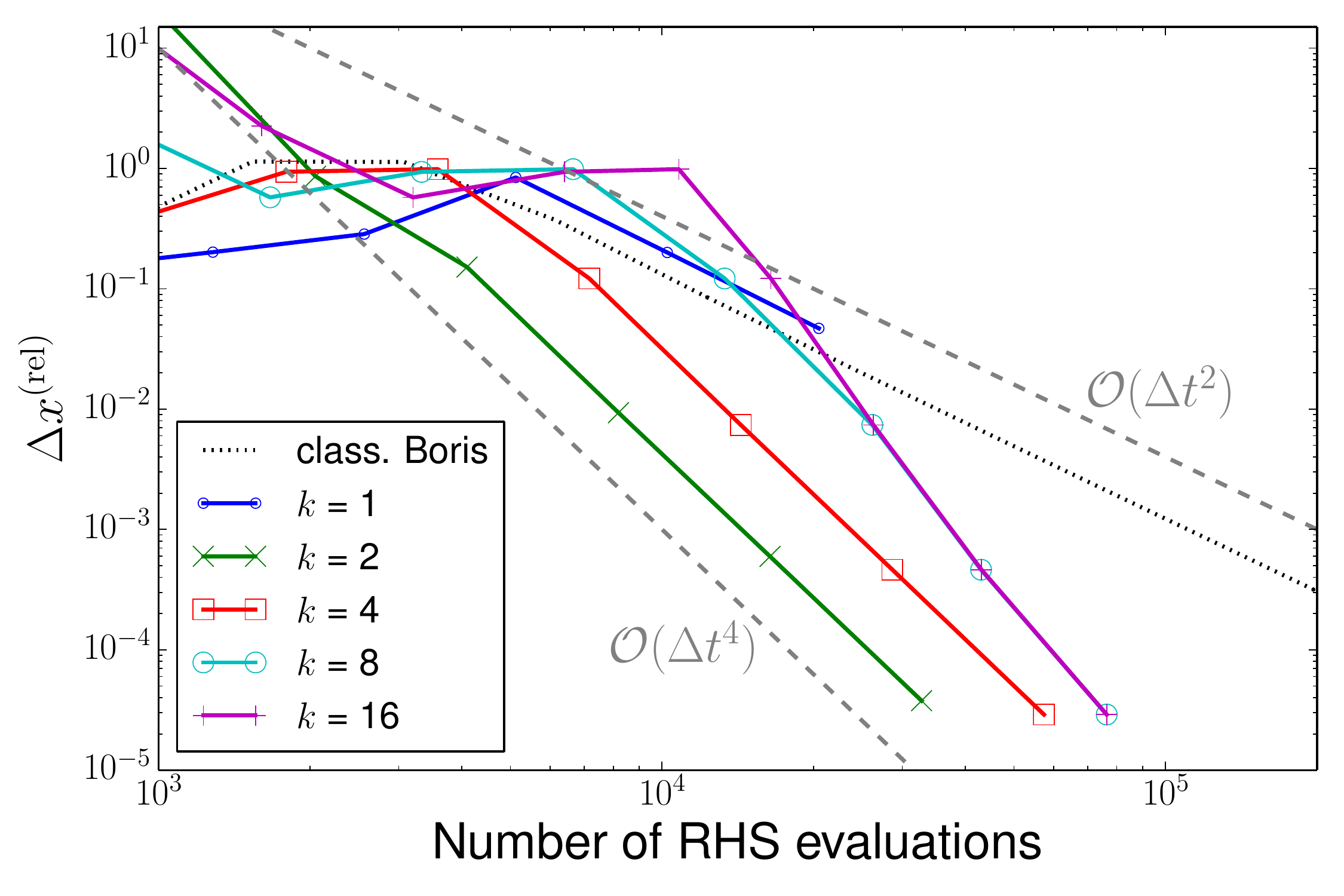}
      \end{minipage}
      \begin{minipage}{\figurewidthTR\textwidth}
          \centerline{$M=5$}
          \includegraphics[width=\textwidth]{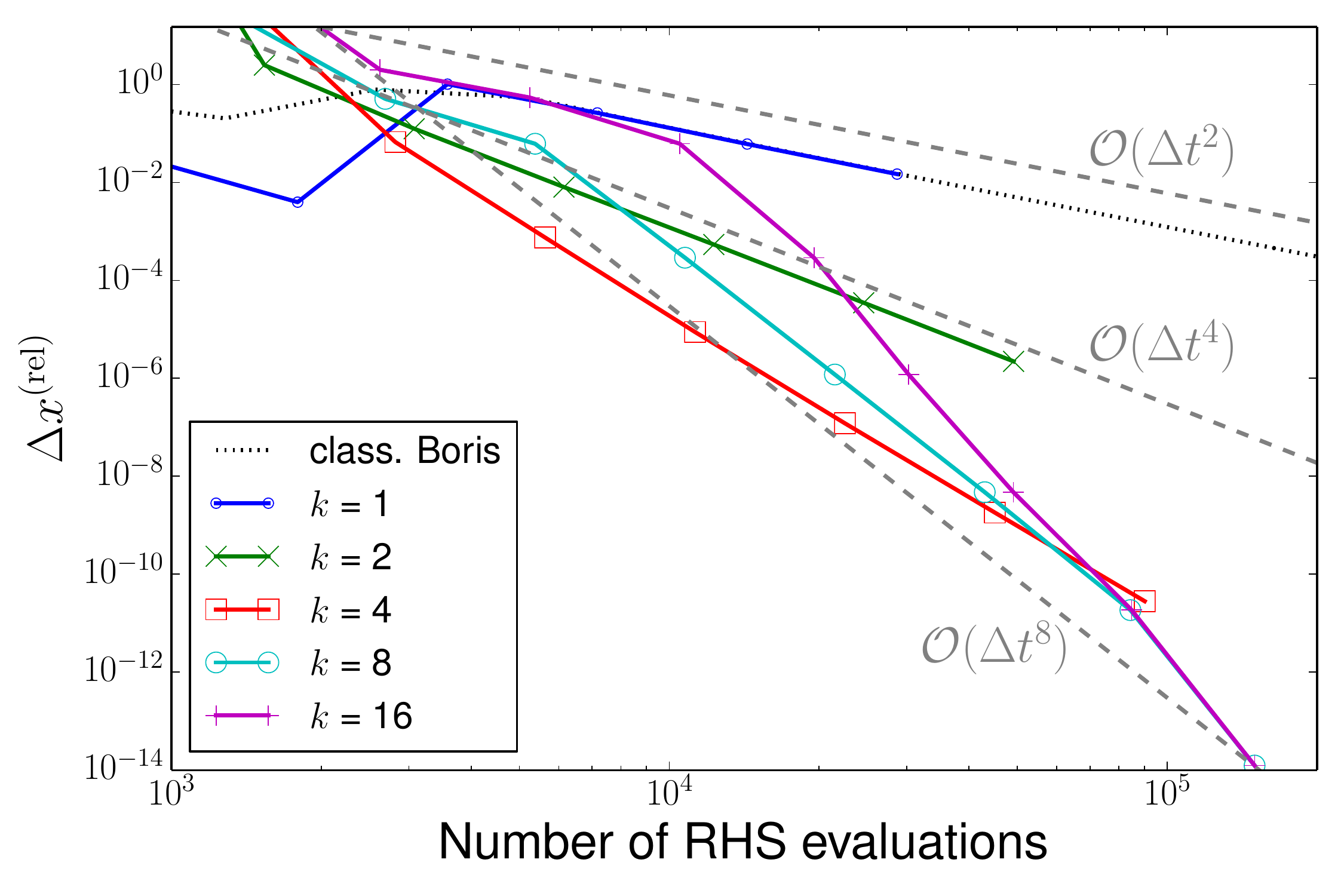}
      \end{minipage}
  \caption{Relative error $\Delta x^{(\text{rel})}$, see~\eqref{eq:rel_error_x}, for the $x$ coordinate of the particle's final position in the Penning trap as a function of the number of right-hand side evaluations performed for 3 and 5 Gauss-Lobatto collocation nodes and different numbers of SDC iterations per time step. The curves for the different runs result from varying the number of time steps for fixed $t_\text{end}$.
  The classical Boris integrator's convergence is shown for comparison.}
  \label{fig:work}
\end{figure}
Boris-SDC requires substantially more work per time step than the original Boris integrator, but can in return provide high-order accuracy.
Figure~\ref{fig:work} shows the error for the classical Boris integrator as well as different configurations of Boris-SDC against the number of evaluations of the right hand side.

For a single iteration and $M=3$, Boris-SDC shows similar precision as the classical Boris solver, because for $k=1$, Boris-SDC corresponds to the classical Boris integrator with smaller time steps, as the time step $[t_n, t_{n+1}]$ is sub-divided by the quadrature nodes.
Increasing the number of right-hand side evaluations by reducing the time step improves accuracy.
For small accuracies, the lower-order methods are the most efficient, that is either the classical Boris integrator or Boris-SDC with one or two iterations only.
For medium to very high accuracy, using more sweeps and thus higher order pays off and requires significantly less right-hand side evaluations than the lower-order versions to reach the same accuracy.

%%%%%%%%%%%%%%%%%%%%%%%%%%%%%%%%%%%%%%%%%%%%%%%%%%%%%%%%%%%%%%%%%%%%%%%%%%%%%%%%
\subsubsection{Energy conservation and symplecticness}\label{sec:energy}
As discussed in Section~\ref{sec:collprop}, the underlying collocation method is symplectic for an appropriate choice of quadrature nodes.
In order to avoid energy drift due the accumulation of round-off errors, see~\cite{Hairer2008}, the computations here have been performed in quadruple precision.

\begin{figure}[t]
      \centering
       \begin{minipage}{\figurewidthTR\textwidth}
          \centerline{$M=3$}
          \includegraphics[width=\textwidth]{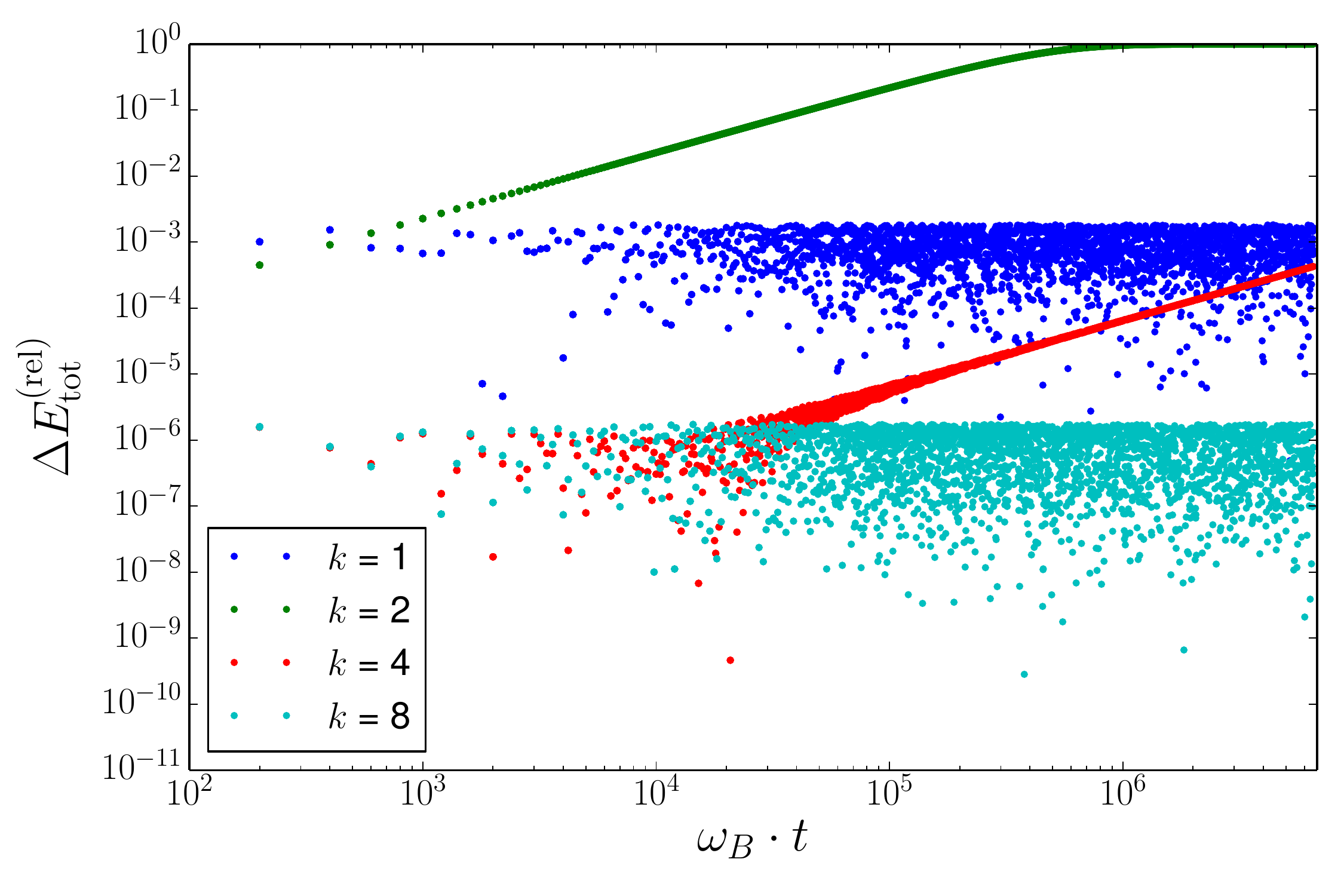}
      \end{minipage}
      \begin{minipage}{\figurewidthTR\textwidth}
          \centerline{$M=5$}
          \includegraphics[width=\textwidth]{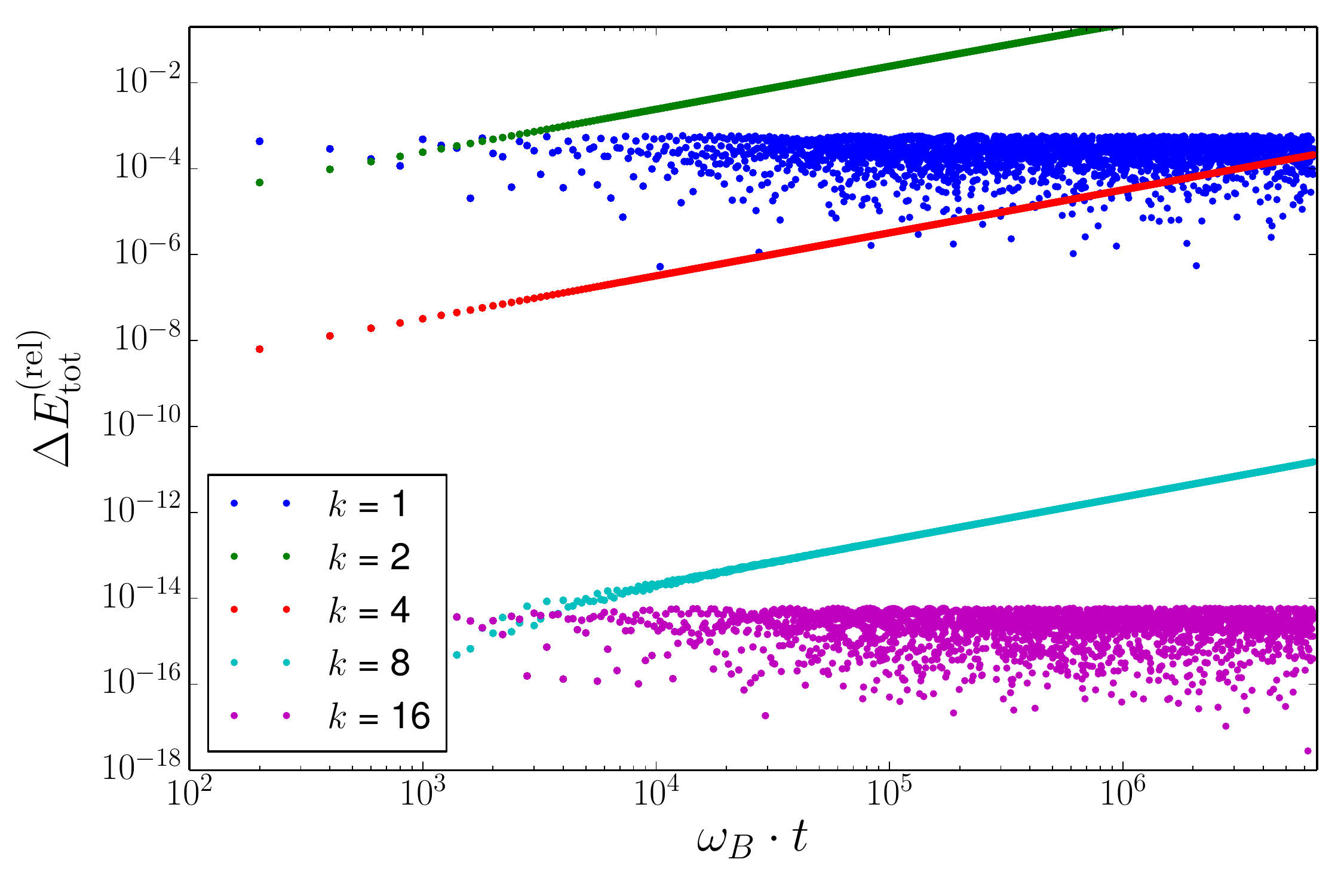}
      \end{minipage}
  \caption{Relative error of the total energy over 16 million time steps for 3 and 5 Gauss-Lobatto collocation nodes and different (but fixed) numbers of SDC iterations per time step.
  All cases were run with quadruple precision to avoid secular drift due to round-off error accumulation, see~\cite{Hairer2008}.
  Since for a single-iteration, Boris-SDC is identical to the classical symplectic velocity-Verlet, the error of the total energy is bounded for long-term simulations there.
  For a small number of iterations but more than one, Boris-SDC has not yet recovered the symplecticness of the underlying collocation method and shows a mild secular drift in the error.
  Once enough SDC iterations have been performed to fully converge to the collocation solution, Boris-SDC retrieves the symplecticness of the collocation method and the energy drift vanishes.
  In addition, due to the high order of Boris-SDC, the energy error is several orders of magnitude smaller than for the classical Boris.}
  \label{fig:energy_gausslobatto}
\end{figure}

Figure~\ref{fig:energy_gausslobatto} shows the relative error in the total energy over $N_\text{steps} = 16'777'216$ time steps ($t_\text{end}=262'144.0$) for $M=3$ and $M=5$ Gauss-Lobatto nodes and different iteration numbers.
All other simulation parameters are identical to those in Table~\ref{tab:param1}.
If only a single sweep is performed, Boris-SDC reduces to the classical Boris-integrator.
Accordingly, the energy error remains bounded over all time steps for both $M=3$ and $M=5$ in this case.
However, because the Boris integrator is only second order accurate, its energy error is comparatively large.
For Boris-SDC,  using more iterations increases the method's accuracy, but at the cost of introducing a slow energy drift:
For small iteration numbers, where Boris-SDC has not yet fully converged, the method is not symplectic and the energy error increases over time.
It grows, however, very slowly:
For four iterations, after 16 million time steps, the energy error is still smaller than for the symplectic second order Boris method for both $M=3$ and $M=5$.
For $M=5$ and eight iterations, although the method is not yet symplectic, the final energy error is still several orders of magnitudes smaller than for the classical Boris.
Moreover, once the number of iterations is set sufficiently large for Boris-SDC to fully converge to the underlying collocation method, symplecticness is retrieved.
For eight iterations for $M=3$ and sixteen iterations for $M=5$,  Boris-SDC no longer shows an energy drift.

  \begin{figure}[t]
      \centering
      \begin{minipage}{\figurewidthWW\textwidth}
          \centerline{$M=3$}
          \includegraphics[width=\textwidth]{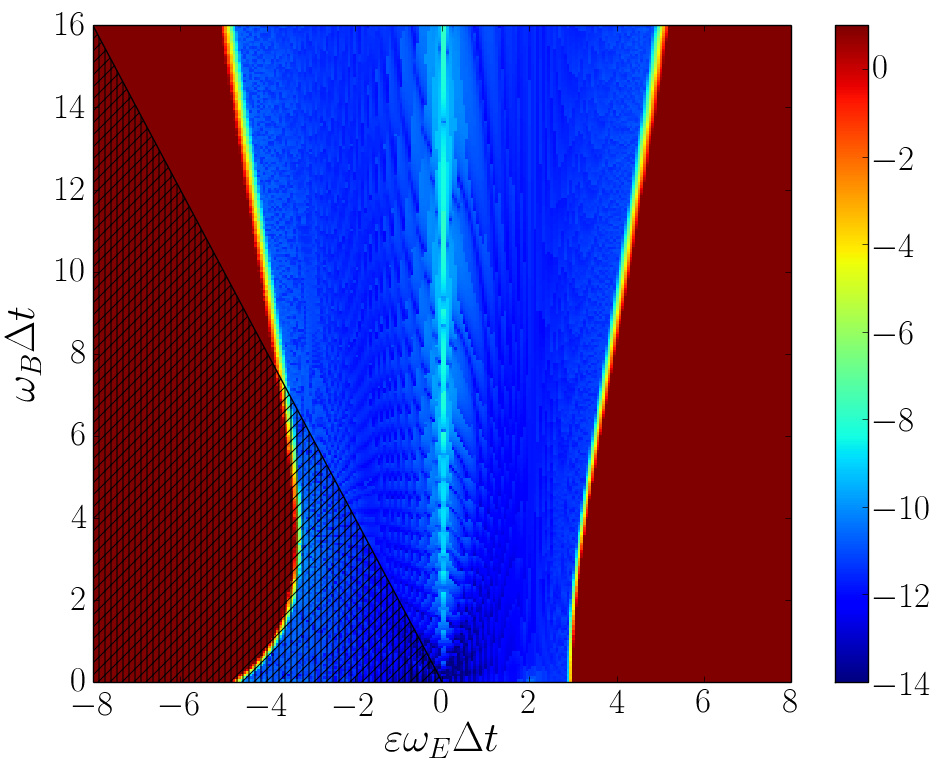}
        \end{minipage}
        \begin{minipage}{\figurewidthWW\textwidth}
            \centerline{$M=5$}
            \includegraphics[width=\textwidth]{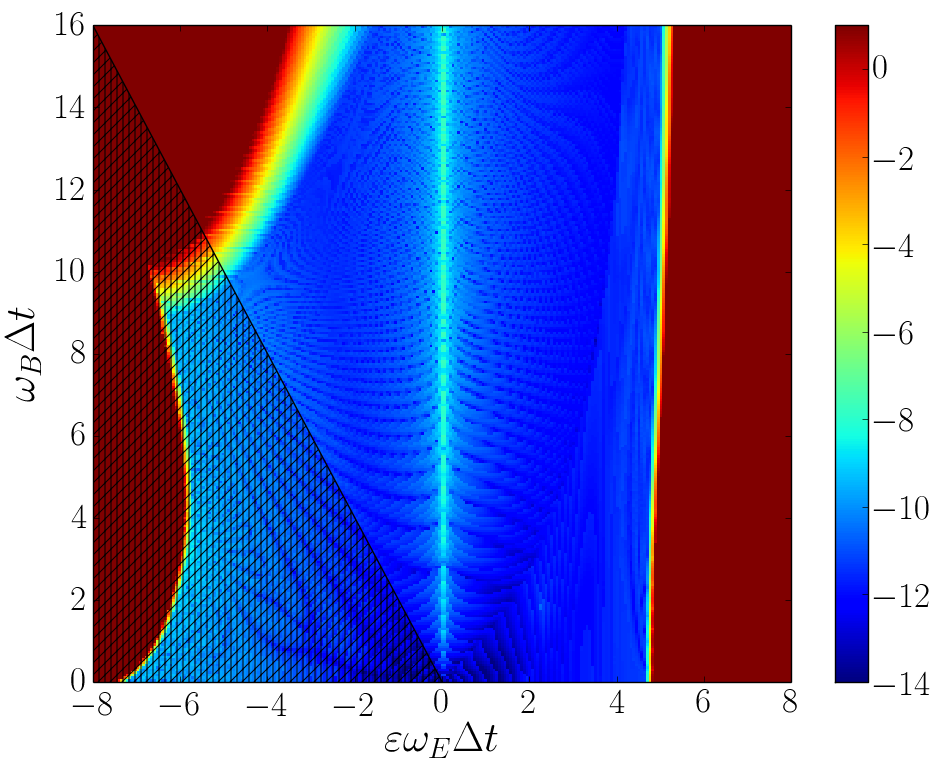}
        \end{minipage}
        \caption{Analysis of the step-to-step energy conservation of Boris-SDC using the update matrix $\tilde{\vect{P}}_\mathrm{sdc}^k$ for $M=3$ and $M=5$ Gauss-Lobatto nodes.
Colors encode $\log_{10}(\max_{i=1\ldots6}|\hat{\vect{H}}_{ii}|)$, compare~\eqref{eq:hhatdef}.
Poor energy conservation is obtained in red, good conservation in blue regions.
These regions essentially correspond to the convergence regions of Figure~\ref{fig:convergence}:
Where Boris-SDC converges, the conservatory properties of the underlying collocation method are reproduced.
The dark gray, hatched region denotes the trap's physical instability region $\omega_B^2 < -4\epsilon\omega_E^2$.
}
        \label{fig:energycons}
    \end{figure}

Instead of fixing the number of iterations, we can, as discussed above, also set a tolerance for the residual and perform sweeps until this tolerance is met.
In this case, a tolerance of $r\leq10^{-10}$ for $M=3$ and $r\leq10^{-18}$ for $M=5$ is required to avoid energy drift.
This illustrates again the point made in~\cite{Hairer2008} that for implicit symplectic methods, accumulating round-off errors in the implicit solver procedure can still lead to energy drift in double precision computations.
The required error tolerance to retrieve symplecticness for $M=5$ here is smaller than machine precision in standard double precision arithmetic.
As completely switching to quadruple precision is computationally expensive on most computers, in~\cite{Hairer2008} other approaches for avoiding this phenomenon are demonstrated.

For the simple linear case studied here, the issue of energy conservation of Boris-SDC can also be studied by analyzing the Boris-SDC update matrix $\tilde{\vect{P}}^{k}_\mathrm{sdc}$.
To this end, we write the system's total energy
\begin{align}
  \mathcal{H} = T + U = \frac{m}{2}\vect{v}^2 + Q\Phi(\vect{x})\label{eq:totenergy}
\end{align}
with particle's mass $m$ and charge $Q$, the kinetic energy $T$, the potential energy $U$ and the external potential $\Phi$ (note that the electric field $\vect{E}=-\nabla\Phi$) as a quadratic form $\mathcal{H} = \vect{u}^T\vect{H}\vect{u}$ with the matrix
\begin{align}
  \vect{H} = \frac{m}{2}
                  \begin{pmatrix}
                     \epsilon\omega_E^2   & 0 & 0 & 0 & 0 & 0 \\
                     0 &  \epsilon\omega_E^2  & 0 & 0 & 0 & 0 \\
                     0 & 0 & -2\epsilon\omega_E^2 & 0 & 0 & 0 \\
                     0 & 0 & 0 & 1 & 0 & 0 \\
                     0 & 0 & 0 & 0 & 1 & 0 \\
                     0 & 0 & 0 & 0 & 0 & 1
                  \end{pmatrix}  
\end{align}
  that acts on the particle's phase space configuration $\vect{u}=(x,y,z,v_x,v_y,v_z)^T$.
  To achieve energy conservation,
  \begin{align}
  \vect{u}_n^T\vect{H}\vect{u}_{n} &= \vect{u}_{n+1}^T\vect{H}\vect{u}_{n+1}\label{eq:energycons}
  \end{align}  
  must hold for all time steps $n$, where the individual steps are linked via the update matrix: $\vect{u}_{n+1} = \tilde{\vect{P}}^{k}_\mathrm{sdc} \vect{u}_n$.
  Note that the iteration count $k$ can be determined via the residual convergence (compare Section~\ref{sec:residual_control}) and is thus in general a function of $n$.
  As in Figure~\ref{fig:convergence}, we use sufficiently many iterations to reach a residual tolerance of $r\leq10^{-12}$.  
  Now,~\eqref{eq:energycons} is fulfilled if
  \begin{align}
    \vect{u}^T\underbrace{\left(\left(\tilde{\vect{P}}_\mathrm{sdc}^k\right)^T\vect{H}\ \tilde{\vect{P}}^{k}_\mathrm{sdc}-\vect{H}\right)}_{\eqcolon\hat{\vect{H}}}\vect{u} &= 0 \label{eq:hhatdef}
  \end{align}
  for all points $\vect{u}$ in phase space that are reachable by regular particle dynamics.
Writing $\vect{u}$ as a linear combination of the standard orthonormal basis in $\mathbb{R}^{6}$ shows that this condition is fulfilled in particular if all diagonal elements of $\hat{\vect{H}}$ are zero.
  
For the same configuration as in Figure~\ref{fig:convergence}, Figure~\ref{fig:energycons} shows $\log_{10}(\max_{i=1\ldots6}|\hat{\vect{H}}_{ii}|)$ for $M=3$ and $M=5$ nodes.
  Small values (blue regions) correspond to small energy error accumulation between time steps while red colors denote bad energy conservation.
The emerging structure is again similar to the convergence regions show in Figure~\ref{fig:convergence}.  
The area of good energy conservation in the parameter space of $\epsilon\omega_E$ and $\omega_B$ is therefore primarily dominated by the convergence properties of the Boris-SDC iteration towards the collocation solution:
  In regions of divergence, energy conservation is violated, while in regions of convergence, the conservatory properties of the underlying collocation are retained.

%%%%%%%%%%%%%%%%%%%%%%%%%%%%%%%%%%%%%%%%%%%%%%%%%%%%%%%%%%%%%%%%%%%%%%%%%%%%%%%%
%%%%%%%%%%%%%%%%%%%%%%%%%%%%%%%%%%%%%%%%%%%%%%%%%%%%%%%%%%%%%%%%%%%%%%%%%%%%%%%%
\subsection{Multiple particles}
 \begin{figure}[t]
      \centering
      \begin{minipage}{\figurewidthTR\textwidth}
          \centerline{$M=3$}
          \includegraphics[width=\textwidth]{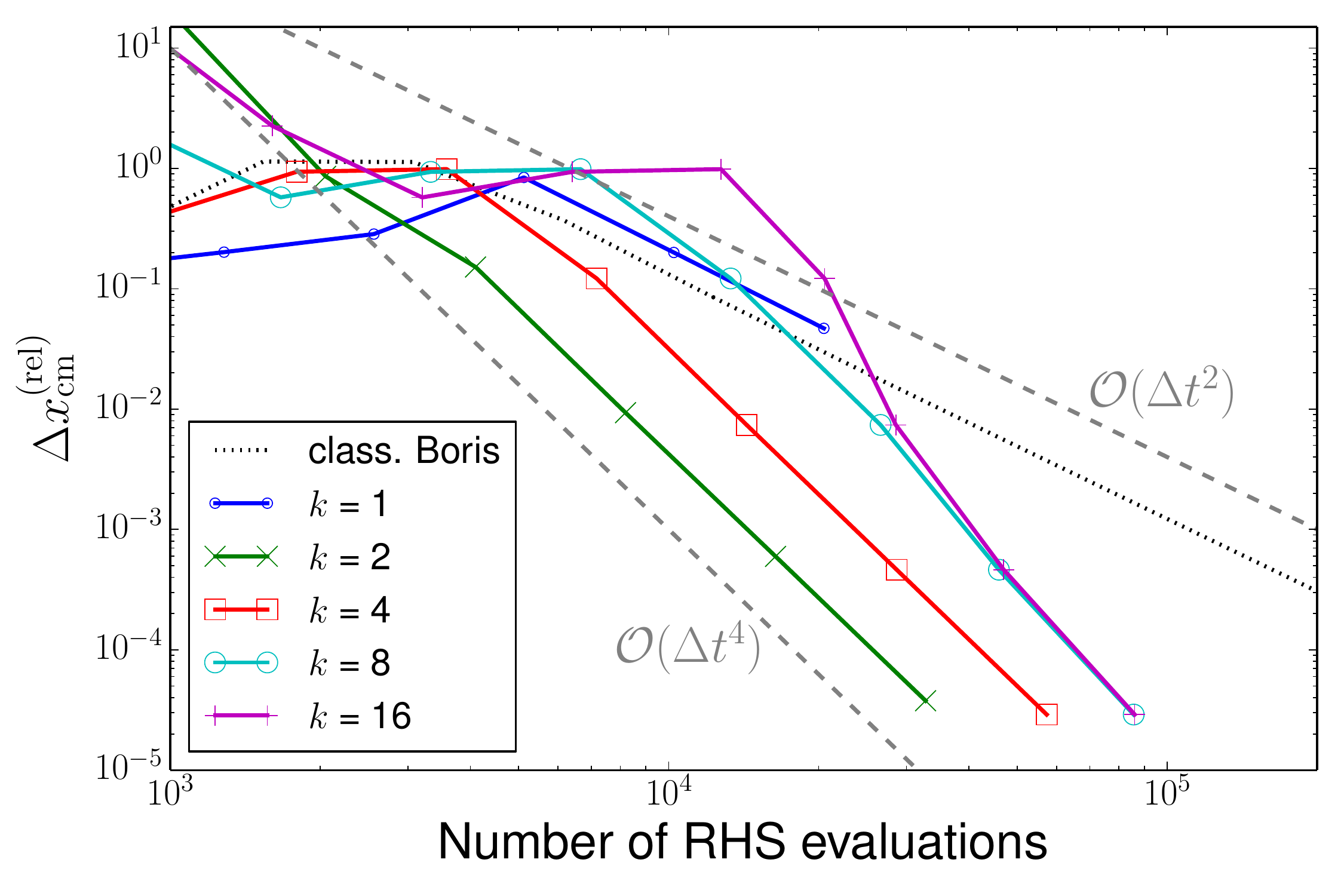}
      \end{minipage}
      \begin{minipage}{\figurewidthTR\textwidth}
          \centerline{$M=5$}
          \includegraphics[width=\textwidth]{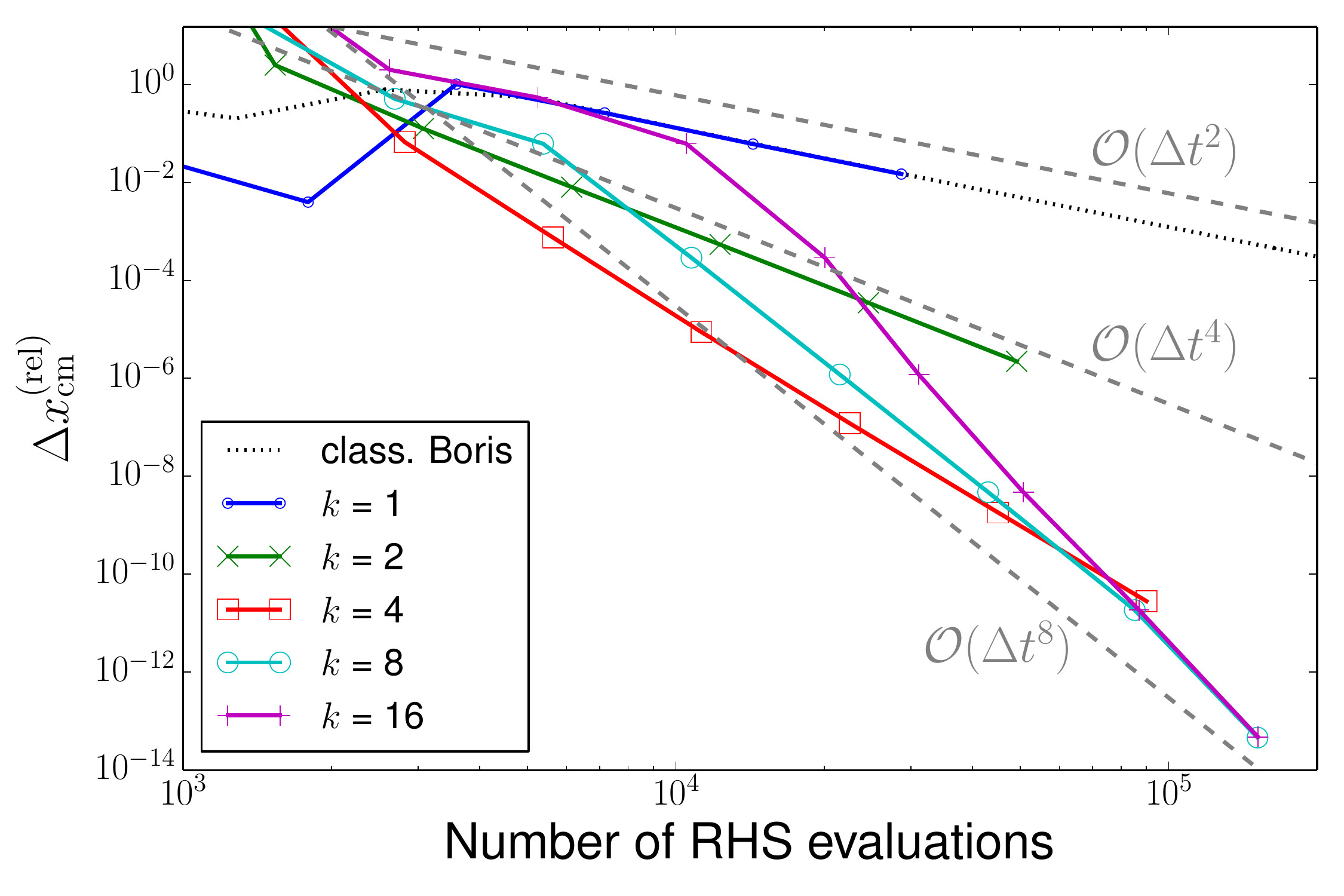}
      \end{minipage}
  \caption{Relative error vs. reference solution for the $x$ coordinate of the particle clouds's centre-of-mass final position in the Penning trap as a function of the number of r.h.s evaluations performed for 3 and 5 Gauss-Lobatto collocation nodes per time step and different number of SDC iterations.
  The curves for the different runs result from varying the total number of time steps for fixed $t_\text{end}$.
  The classical Boris integrator's convergence is shown for comparison.
  }
  \label{fig:work100iter}
\end{figure}
In this example, instead of a single particle we study a cloud of $100$ particles inside the trap.
The setup parameters are identical to those in Table~\ref{tab:param1}.
However, the individual particle's initial positions and velocities are distorted by random vectors $|\vec{x}_\text{shift}|\leq0.001$ and $|\vec{v}_\text{shift}|\leq5.0$, respectively.
To ensure comparability, the same random distortions are applied across runs with different methods.
Because no analytical solution is available for the particle cloud, a reference solution is computed with a high-order run using a very fine time step. % reference data was computed with Boris-SDC, $M=9$ nodes, $N_\mathrm{steps}=131072$, 16 iterations

Instead of the error in position and velocity of individual particles, we track the relative error of the position of the center-of-mass
\begin{align}
  \vect{x}_\text{cm} &= \frac{\sum_{i=1}^{N_\text{particles}} m_i \vect{x}_i}{\sum_{i=1}^{N_\text{particles}} m_i}
\end{align}
of the particle cloud.
Figure~\ref{fig:work100iter} shows this error versus the number of right-hand side evaluations for Boris-SDC with different numbers of iterations.
As in the single particle case, each sweep increases the order of convergence by about two, although for $N=5$ with four sweeps, the result is rather seventh than eighth order accurate.
The higher order of Boris-SDC also pays off for the particle cloud:
For a medium to high precision simulation of the center-of-mass, it requires significantly fewer evaluations of the right hand side than the classical Boris integrator to achieve the same accuracy.

 \begin{figure}[t]
      \centering
      \begin{minipage}{\figurewidthTR\textwidth}
          \centerline{$M=3$}
          \includegraphics[width=\textwidth]{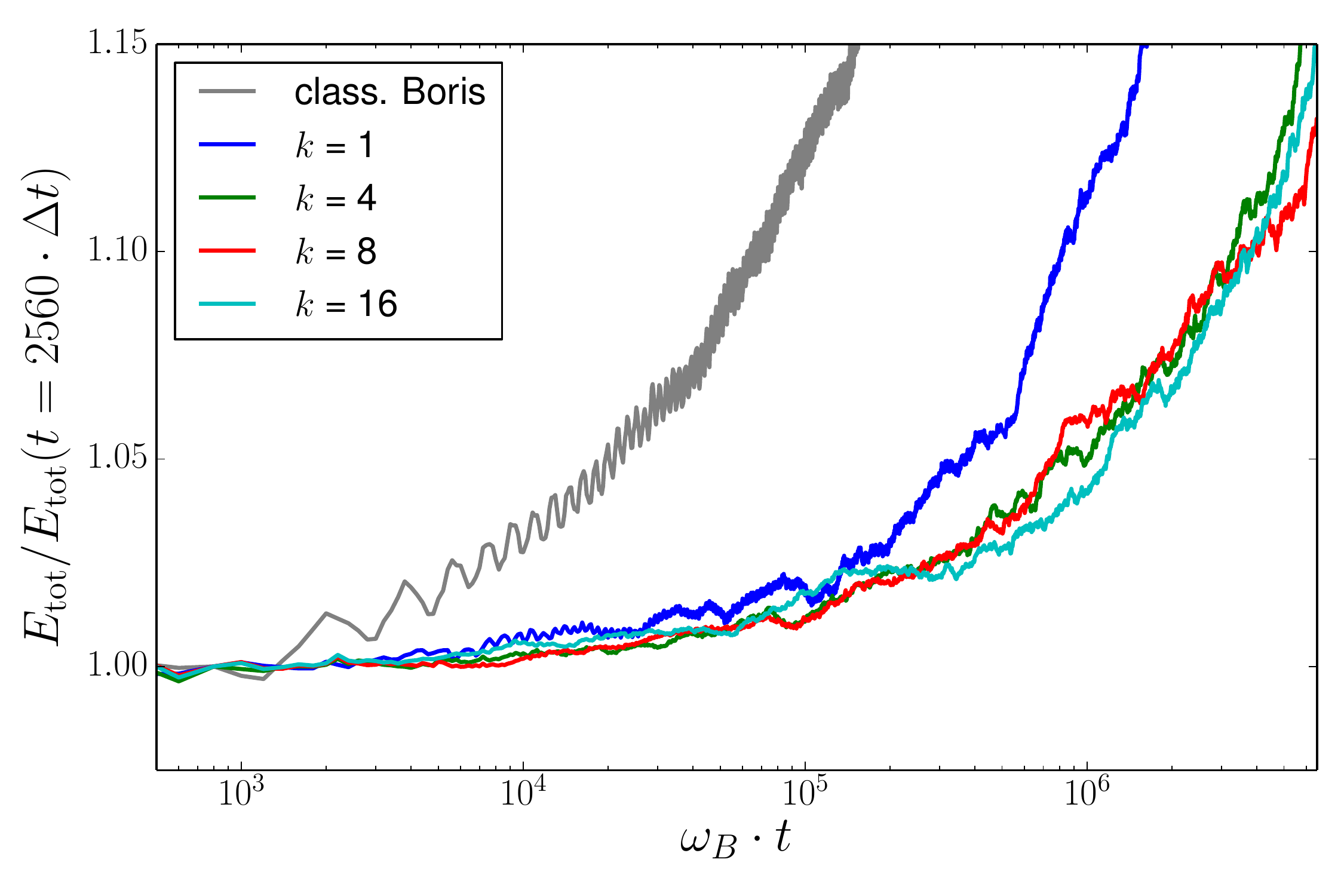}
      \end{minipage}
      \begin{minipage}{\figurewidthTR\textwidth}
          \centerline{$M=5$}
          \includegraphics[width=\textwidth]{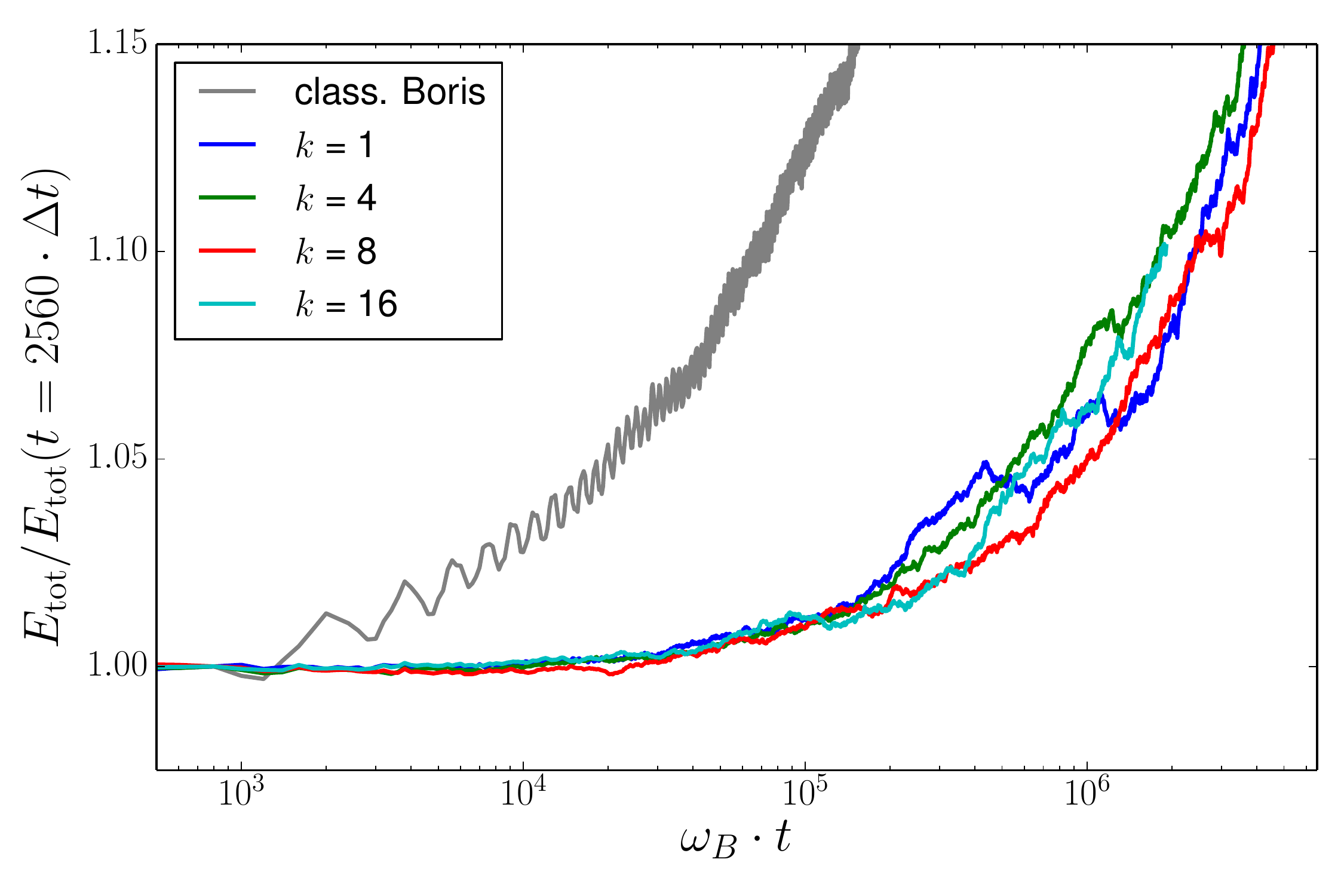}
      \end{minipage}
  \caption{Study on long-term energy stability for the multi-particle setup.
           To avoid the effects of the initial non-equilibrium (random) setup, the ratio of the total energy to a reference value after an initial relaxation phase of 2560 time steps is shown.
           Over $N_\text{steps}=16'777'216$ simulation steps ($t_\text{end}=262'144.0$) a noticeable energy drift is observed due to under-resolving close encounters between the simulated particles.
           This effect is significantly less pronounced for Boris-SDC.
          }
  \label{fig:energy100}
\end{figure}
For a molecular dynamics simulations with Coulomb interaction, significant numerical heating occurs in the beginning due to close encounters of particles.
To assess energy conservation, we therefore compute in every run first an initial relaxation phase of 2'560 time steps and use the energy at the end of the initial phase as reference.
Figure~\ref{fig:energy100} shows the ratio of the systems total energy~\eqref{eq:totenergy} to this reference value for the same simulation time step and total runtime as used for Figure~\ref{fig:energy_gausslobatto}.
Because of the rather coarse time step, the original Boris method exhibits a significant energy drift from under-resolved Coulomb collisions.
The increased complexity of Boris-SDC, which evaluates forces on particles at several intermediate steps, results in improved energy stability due to improved trajectory resolution:
The energy drift for Boris-SDC sets in much later than for the classical Boris integrator, leading to a reduction of numerical heating.

\section{Conclusions and Outlook}\label{sec:conclusion}

The Boris integration method is the de-facto standard approach for computing trajectories of particles in electric and magnetic fields and is widely used in a large variety of applications.
While it is easy to implement and lightweight in computational cost, it is only second-order accurate in time and an extension to higher orders is not straightforward.

In this work, we introduced Boris-SDC, which combines the classical Boris method with a flexible, iterative, spectral deferred corrections solver for the collocation formulation of second-order ODEs.
The derivation extended the SDC approach for Newtonian equations of motion to the case of velocity-dependent forces.
It also showed that the arising implicit system can be cast into the same form as in the Boris method and therefore be solved in an analogue way.
Boris-SDC is thus specifically tailored for simulations of particles in electric and magnetic fields; it maintains the advantage of the Boris integrator to be essentially explicit while allowing to generically construct a method of arbitrary order by varying the number of SDC iterations and the number of quadrature nodes.
Whether Boris-SDC performs better than fully implicit methods remains to be seen and is not a straightforward question to answer: This will depend on the chosen problem and on the used nonlinear and linear solver for the implicit method.
Corresponding comparisons are left for future work.

The properties of Boris-SDC were studied through numerical examples and compared to the classical Boris method for different particle setups in a classical Penning trap.
It was shown that Boris-SDC has a larger stability region and provides high-order accuracy for a single particle trajectory and the center-of-mass of a particle cloud.
In particular, we observe that in general each SDC iteration with the second-order velocity-Verlet method also increases the order of Boris-SDC by two.
For medium to high accuracies, Boris-SDC requires significantly fewer force evaluations because of the higher order.
For sufficiently many iterations, Boris-SDC replicates the excellent conservation properties of the underlying collocation method.
We also demonstrated that the properties of Boris-SDC as a time stepping method are intimately linked to the convergence behavior of the SDC iteration as a preconditioned Picard iteration to solve the collocation equation.

This work provides an illustrative example how the interpretation of SDC as preconditioned Picard iteration can be used to incorporate new and rather complex base methods which exploit specific features of a certain problem.
As a next step, Boris-SDC can be incorporated and tested in a legacy code for a realistic, real-world application.
The fields of high-intensity laser plasma interaction, particle trajectory integration in accelerators and space-weather studies are generic candidates.
The non-intrusive nature of our approach facilitates the augmentation of an existing Boris integrator, e.\,g.~in a legacy Particle-In-Cell code, with Boris-SDC.
A comparison of Boris-SDC with other high order methods would be an interesting direction of future research as well.

Furthermore, the Boris-SDC method we introduced can provide the foundation for the development of a new variant of the time-parallel PFASST method (see~\cite{EmmettMinion2012,EmmettMinion2014_DDM}), optimized for particle simulations in plasma physics.
Thus, the integration of Boris-SDC into the novel \texttt{PFASST++} framework~\cite{PFASST++} is planned for future work.
This requires extending Boris-SDC to multiple levels in space and time with adequate coarsening strategies, see~\cite{SpeckEtAl2014_BIT} for a description of multi-level SDC and~\cite{SpeckEtAl2012,SpeckEtAl2014_DDM2012} for a first idea for particle-based coarsening.
An efficient time-parallel method tailored for particle simulations could greatly aid in better exploiting the computational resources of massively parallel high-performance computing systems for plasma physics applications.
Extending Boris-SDC to Boris-PFASST could provide the keystone for large-scale space-time parallel particle simulations in plasma physics.

\appendix

\section{Velocity-Verlet integration in matrix formulation}\label{sec:apx_vv}

The notation discussed here is based on the ideas presented in~\cite{Minion2014_2ndOrderSDC}. 
In addition, we make use of the definitions introduced in Section~\ref{sec:coll}.
We consider a particle with position $\vect{x}\in\mathbb{R}^d$ and velocity $\vect{v}\in\mathbb{R}^d$ at time $t>0$.
Note, that the dimension of space $d$ is just a formal notation.
It can also be extended to represent the coordinates of a multitude of individual particles.
Newton's equations of motion for the particle are then given by
\begin{align}
	\frac{\d\vect{v}}{\d t} = \vect{f}(\vect{x},\vect{v}),\quad \frac{\d\vect{x}}{\d t} = \vect{v}\label{eq:apx_eom}
\end{align}
with a suitable right-hand side $\vect{f}\in\mathbb{R}^d$.
Using standard velocity-Verlet integration with time steps $\tau_0,\ldots,\tau_M$, $M\ge1$, with $\Delta\tau_{m+1} = \tau_{m+1}-\tau_{m}$, $m=0,\ldots,M-1$, this can be discretized as
\begin{subequations}
\label{eq:apx_vv}
\begin{align}
	\vect{v}_{m+1} &= \vect{v}_m + \frac{\Delta\tau_{m+1}}{2}\left(\vect{f}_m + \vect{f}_{m+1}\right)\label{eq:apx_vv_v}\\
	\vect{x}_{m+1} &= \vect{x}_m + \Delta\tau_{m+1}\left(\vect{v}_m + \frac{\Delta\tau_{m+1}}{2}\vect{f}_m\right)\label{eq:apx_vv_x}
\end{align}
\end{subequations}
where $\vect{x}_m \approx \vect{x}(\tau_m)$, $\vect{v}_m \approx \vect{v}(t_m)$ and $\vect{f}_m = \vect{f}(\vect{x}_m,\vect{v}_m)$. 
The update formula for $\vect{x}$ consists of an Euler half-step for $\vect{v}$ plus an Euler full-step for $\vect{x}$ and is explicit. 
The update formula for $\vect{v}$, however, is given by a trapezoidal rule, which is implicit due to the $\vect{v}$-dependency of $\vect{f}$ and second-order accurate.

To describe the process of velocity-Verlet integration from $\tau_0$ to $\tau_M$, we take Eqs.~\eqref{eq:apx_vv} and apply these formulas recursively, obtaining
\begin{align}
	\vect{v}_{m+1} &= \vect{v}_0 + \frac{1}{2}\sum_{l=1}^{m+1}\Delta\tau_l\left(\vect{f}_{l-1} + \vect{f}_l\right)\label{eq:apx_vv_v0}\\
	\vect{x}_{m+1} &= \vect{x}_0 + \sum_{l=1}^{m+1}\Delta\tau_l\vect{v}_{l-1} + \frac{1}{2}\sum_{l=1}^{m+1}(\Delta\tau_l)^2\vect{f}_{l-1}. \label{eq:apx_vv_x0}
\end{align}
We now introduce matrices
\begin{align}
	Q_E &\coloneq 
	\begin{pmatrix}
	0 & 0 & 0 & \cdots & 0\\
	\Delta\tau_1 & 0 & 0 \\
	\Delta\tau_1 & \Delta\tau_2 & 0 & & \vdots\\
	\vdots & \vdots & \ddots & \ddots \\
	\Delta\tau_1 & \Delta\tau_2 & \cdots & \Delta\tau_{M} & 0
	\end{pmatrix},&
	Q_I &\coloneq 
	\begin{pmatrix}
	0 & 0 & 0 & \cdots & 0\\
	0 & \Delta\tau_1 & 0  \\
	0 & \Delta\tau_1 & \Delta\tau_2 & & \vdots\\
	\vdots & \vdots & \ddots & \ddots \\
	0 & \Delta\tau_1 & \Delta\tau_2 & \cdots & \Delta\tau_{M} 
	\end{pmatrix}
\end{align}
and 
\begin{align}
	Q_T \coloneq \frac{1}{2}(Q_E + Q_I),\label{eq:QT}
\end{align}
representing the propagation matrices for the explicit Euler, the implicit Euler and the trapezoidal rule.
We have $Q_E,Q_I,Q_T\in\mathbb{R}^{(M+1)\times(M+1)}$.
Furthermore, we gather the values $\vect{v}_m,\vect{x}_m,\vect{f}(\vect{x}_m,\vect{v}_m)\in\mathbb{R}^d$ into vectors $\tvect{v}, \tvect{x}, \tvect{f}(\tvect{X},\tvect{V})$ as in Section~\ref{sec:coll}.
Then, Eqs.~\eqref{eq:apx_vv_v0} and~\eqref{eq:apx_vv_x0} with the notation defined in Section~\ref{sec:coll} read
\begin{subequations}
\label{eq:apx_vv_vec}
\begin{align}
	\tvect{v} &= \tvect{v}_0 + \vect{Q}_T\tvect{f}(\tvect{X},\tvect{V})\label{eq:apx_vv_v0_vec}\\
	\tvect{x} &= \tvect{x}_0 + \vect{Q}_E\tvect{v} + \frac{1}{2}(\vect{Q}_E\circ\vect{Q}_E)\tvect{f}(\tvect{X},\tvect{V})\label{eq:apx_vv_x0_vec}
\end{align}
\end{subequations}
with the Hadamard product $\circ$ (entry-wise multiplication of two matrices) and identity matrix $\vect{I}_d\in\mathbb{R}^{d\times d}$.
%Note that using this notation the initial values $\vect{v}_0,\vect{x}_0$ are always retained due to the first row of zeros in $\vect{Q}_E$ and $\vect{Q}_I$.
We can now use~\eqref{eq:apx_vv_v0_vec} to rewrite~\eqref{eq:apx_vv_x0_vec} as
\begin{align}
	\tvect{x} = \tvect{x}_0 + \vect{Q}_E\tvect{v}_0 + \vect{Q}_x\tvect{f}(\tvect{X},\tvect{V})\label{eq:apx_vv_x0v0_vec}
\end{align}
with matrix
\begin{align}
	\vect{Q}_x \coloneq \vect{Q}_E\vect{Q}_T+\frac{1}{2}(\vect{Q}_E\circ\vect{Q}_E).\label{eq:Qx}
\end{align}
With~\eqref{eq:apx_vv_v0_vec} and~\eqref{eq:apx_vv_x0v0_vec}, both $\tvect{v}$ and $\tvect{x}$ solely depend on the initial vectors $\tvect{v}_0$ and $\tvect{x}_0$ as well as the function vector $\tvect{f}(\tvect{X},\tvect{V})$.

In order to combine both equations into a single formula based on $\tvect{u}$ (see Section~\ref{sec:coll}), we make again use of the permutation operators $\krI{I}{x}$, $\krI{I}{v}$ and $\krI{I}{xv}$ of \eqref{eq:I_matrices} and define
\begin{align}
	\vect{Q}_\mathrm{vv} &\coloneq \vect{\krI{Q_{\mathnormal x}}{x}} + \vect{\krI{Q_{\mathnormal T}}{v}} = \left(Q_EQ_T+\frac{1}{2}(Q_E\circ Q_E)\right)\otimes \krI{I}{x}\otimes\vect{I}_d + Q_T\otimes \krI{I}{v}\otimes \vect{I}_d.
\end{align}
Then, the combined formulation of the velocity-Verlet scheme for $M$ substeps in matrix formulation is given by
\begin{align}
	\tvect{u} =\ \tvect{u}_0 + \vect{\krI{Q_{\mathnormal E}}{xv}}\tvect{u}_0 + \vect{Q}_\mathrm{vv}\tvect{f}(\tvect{u}) \coloneq
	\vect{C}_\mathrm{vv}\tvect{u}_0 + \vect{Q}_\mathrm{vv}\tvect{f}(\tvect{u})\label{eq:apx_vv_mat}
\end{align}
with $\vect{C}_\mathrm{vv} \coloneq \vect{I}_{(M+1)2d} + \vect{\krI{Q_{\mathnormal E}}{xv}}$.
Formally, we write this as (probably non-linear) system
\begin{align}
	\vect{M}_\mathrm{vv}(\tvect{u}) = \vect{C}_\mathrm{vv}\tvect{u}_0\label{eq:apx_vv_system}
\end{align}
with
\begin{align}
	\vect{M}_\mathrm{vv}(\cdot) \coloneq \left(\vect{I}_{(M+1)2d} -  \vect{Q}_\mathrm{vv}\tvect{f}\right)(\cdot),\label{eq:apx_vv_systemoperator}
\end{align}
so that the solution vector $\tvect{u}$ is formally given by $\tvect{u} = \inv{\vect{M}}_\mathrm{vv}(\vect{C}_\mathrm{vv}\tvect{u}_0)$.
Using the linear transfer operators $\vect{T}_\mathrm{P}\in\mathbb{R}^{(M+1)2d\times 2d}$ and $\vect{T}_\mathrm{R}\in\mathbb{R}^{2d\times (M+1)2d}$ of~\eqref{eq:transfer_ops}, we finally obtain the formal update formulation for the velocity-Verlet scheme over $M$ substeps with
\begin{align}
	\vect{u}_M = \vect{P}_\mathrm{vv}(\vect{u}_0) \coloneq \vect{T}_\mathrm{R} \inv{\vect{M}}_\mathrm{vv}(\vect{C}_\mathrm{vv}\vect{T}_\mathrm{P}\ \vect{u}_0).\label{eq:apx_vv_update}
\end{align}
%We note that for the special case of linear right-hand sides $\vect{f}(\vect{x},\vect{v}) = \vect{A}(\vect{x},\vect{v})\in\mathbb{R}^d$, i.\,e.~with matrix $\vect{A}\in\mathbb{R}^{d\times 2d}$, the system operator becomes a matrix with
%\begin{align} 
%	\vect{M}_\mathrm{vv} = \vect{I}_{(M+1)2d} -  (\vect{Q}_\mathrm{vv}\otimes\vect{A})
%\end{align}
%and the update formula reads
%\begin{align}
%	\vect{u}_M = \vect{P}_\mathrm{vv}\vect{u}_0 \coloneq \vect{T}_\mathrm{R} \inv{\vect{M}}_\mathrm{vv}  \vect{C}_\mathrm{vv}\vect{T}_\mathrm{P}\ \vect{u}_0.
%\end{align}
Applying $\inv{\vect{M}}_\mathrm{vv}$ corresponds to solving the non-linear system~\eqref{eq:apx_vv_system}.
The special lower block-diagonal structure of $\vect{M}_\mathrm{vv}$ (containing diagonal elements due to the implicit dependency in $\vect{Q}_T$) makes inversion easy:
As for the classical velocity-Verlet notation, each solution $\vect{u}_m$ is obtained from the previous ones. 
Thus, applying $\inv{\vect{M}}_\mathrm{vv}$ simply means stepping from $\tau_0$ to $\tau_M$ using standard velocity-Verlet integration.

%%%%%%%%%%%%%%%%%%%%%%%%%%%%%%%%%%%%%%%%%%%%%%%%%%%%%%%%%%%%%%%%%%%%%%%%%%%%%%%%
%%%%%%%%%%%%%%%%%%%%%%%%%%%%%%%%%%%%%%%%%%%%%%%%%%%%%%%%%%%%%%%%%%%%%%%%%%%%%%%%
%%%%%%%%%%%%%%%%%%%%%%%%%%%%%%%%%%%%%%%%%%%%%%%%%%%%%%%%%%%%%%%%%%%%%%%%%%%%%%%%
%%%%%%%%%%%%%%%%%%%%%%%%%%%%%%%%%%%%%%%%%%%%%%%%%%%%%%%%%%%%%%%%%%%%%%%%%%%%%%%%
\section*{Acknowledgments}
This work greatly benefited from the \texttt{libpfasst} library developed by M.~Minion and M.~Emmett which provided the infrastructure for the implementation of Boris-SDC. 
We also gratefully acknowledge many inspiring discussions with both colleagues as well as with the members of the Simulation Laboratory Plasma Physics, in particular the \texttt{pepc} developer's group at J\"ulich Supercomputing Centre.
Numerical experiments were performed on the JUROPA (through computing time grant~JZAM04) and JUDGE systems in J\"ulich.
Robert Speck and Daniel Ruprecht acknowledge support by Swiss National Science Foundation grant 145271 under the lead agency agreement through the project "ExaSolvers" within the Priority Programme 1648 "Software for Exascale Computing" of the Deutsche Forschungsgemeinschaft.

\bibliography{Pint,Pint_Self,sdc,refs,mwlibrary}

\end{document}